\documentclass[11pt]{amsart}
\usepackage{amsmath,amsthm, amscd, amssymb, amsfonts}
\usepackage[all]{xy}
\usepackage{hhline}

\usepackage{ifthen}

\newcommand{\Ind}{\operatorname{Ind}}
\newcommand{\ord}{\operatorname{ord}}

\newcommand\toba{{\mathfrak B }}
\newcommand\trasp{\pi}

\newcommand{\trid}{\triangleright}

\newcommand{\ku}{\mathbb C}

\newcommand{\Z}{{\mathbb Z}}
\newcommand{\N}{{\mathbb N}}

\newcommand{\G}{{\mathbb G}}

\newcommand{\C}{{\mathcal C}}

\newcommand{\B}{{\mathcal B}}

\newcommand{\Oc}{{\mathcal O}}

\newcommand\card{\operatorname{card}}

\newcommand\sgn{\operatorname{sgn}}

\theoremstyle{plain}

\newtheorem{maintheorem}{Theorem}
\newtheorem{lema}{Lemma}[section]
\newtheorem{theorem}[lema]{Theorem}

\newtheorem{prop}[lema]{Proposition}

\theoremstyle{definition}
\newtheorem{definition}[lema]{Definition}

\theoremstyle{remark}
\newtheorem{obs}[lema]{Remark}

\newcommand\id{\operatorname{id}}

\newcommand\Id{\operatorname{Id}}

\newcommand\sn{\mathbb S_n}
\newcommand\snn{\mathbb S_{2n}}

\newcommand{\Dtriangle}[7]{
\rule[-3\unitlength]{0pt}{12\unitlength}
\begin{picture}(18,7)(0,3)
\put(4,4){\ifthenelse{\equal{#1}{l}}{\circle*{2}}{\circle{2}}}
\put(5,4){\line(1,0){8}}
\put(14,4){\ifthenelse{\equal{#1}{r}}{\circle*{2}}{\circle{2}}}
\put(4.4472,4.4472){\line(1,2){4.1056}}
\put(9,14){\ifthenelse{\equal{#1}{t}}{\circle*{2}}{\circle{2}}}
\put(13.5528,4.8944){\line(-1,2){4.1056}}
\put(2,3){\makebox[0pt][r]{\scriptsize #2}}
\put(9,17){\makebox[0pt]{\scriptsize #3}}
\put(16,3){\makebox[0pt][l]{\scriptsize #4}}
\put(6,9){\makebox[0pt][r]{\scriptsize #5}}
\put(12.5,9){\makebox[0pt][l]{\scriptsize #6}}
\put(9,1){\makebox[0pt]{\scriptsize #7}}
\end{picture}}
\newcommand{\cDtriangle}[7]{
\rule[-10\unitlength]{0pt}{\unitlength}
\begin{picture}(18,7)(0,10)
\put(4,4){\ifthenelse{\equal{#1}{l}}{\circle*{2}}{\circle{2}}}
\put(5,4){\line(1,0){8}}
\put(14,4){\ifthenelse{\equal{#1}{r}}{\circle*{2}}{\circle{2}}}
\put(4.4472,4.8944){\line(1,2){4.1056}}
\put(9,14){\ifthenelse{\equal{#1}{t}}{\circle*{2}}{\circle{2}}}
\put(13.5528,4.8944){\line(-1,2){4.1056}}
\put(2,3){\makebox[0pt][r]{\scriptsize #2}}
\put(9,17){\makebox[0pt]{\scriptsize #3}}
\put(16,3){\makebox[0pt][l]{\scriptsize #4}}
\put(6,9){\makebox[0pt][r]{\scriptsize #5}}
\put(12.5,9){\makebox[0pt][l]{\scriptsize #6}}
\put(9,1){\makebox[0pt]{\scriptsize #7}}
\end{picture}}

\def\pf{\begin{proof}}
\def\epf{\end{proof}}

\theoremstyle{remark}

\begin{document}

\renewcommand{\baselinestretch}{1.2}

\thispagestyle{empty}
%\vspace*{2in}

\title[pointed Hopf algebras with coradical $\ku\sn$]
{On pointed Hopf algebras associated to unmixed conjugacy classes
in $\sn$}
\author[andruskiewitsch and fantino]{Nicol\'{a}s Andruskiewitsch}
\address{Facultad de Matem\'{a}tica, Astronom\'{i}a y F\'{i}sica\\
Universidad Nacional de C\'{o}rdoba \\ CIEM - CONICET,
(5000) Ciudad Universitaria \\
C\'{o}rdoba \\Argentina}\email{andrus@famaf.unc.edu.ar}
\author[]{Fernando Fantino}
\email{fantino@famaf.unc.edu.ar}
\thanks{This work was partially
supported by CONICET, ANPCyT and Secyt (UNC)}

\date{August 28, 2006}

\begin{abstract}
Let $\pi \in \sn$ be a product of disjoint cycles of the same
length, $\C$ the conjugacy class of $\pi$ and $\rho$ an
irreducible representation of the isotropy group of $\pi$. We
prove that either the Nichols algebra $\toba(\C,\rho)$ is
infinite-dimensional, or the braiding of the Yetter-Drinfeld
module is negative.
\end{abstract}
\maketitle

\setcounter{tocdepth}{2}

\tableofcontents

\section*{Introduction}\label {0}

Hopf algebras have important applications in mathematics and
mathematical physics. Indeed, Hopf algebras give rise to finite
tensor categories in the sense of \cite{ENO, EO} through their
categories of representations. In this way, for instance,
semisimple Hopf algebras are present in a fundamental way in
rational conformal field theories. Also, non-semisimple Hopf
algebras are related to logarithmic conformal field theories
\cite{Ga}. It is natural to expect that classification results on
finite-dimensional Hopf algebras would have a significant impact
in those areas. Needless to say, classification efforts often come
out in discovery of new examples.

This article, in the line of \cite{AZ}, is a contribution to the
classification of finite-dimensional complex pointed Hopf algebras
$H$ with $G(H)$ non-abelian, in the framework of the Lifting
Method \cite{AS-cambr}. Let $G$ be a finite group. An important
stage in the proposed approach to classification of
finite-dimensional complex pointed Hopf algebras $H$ with $G(H) =
G$ is the determination of all Yetter-Drinfeld modules $V$ over
the group algebra of $G$ such that the Nichols algebra $\toba(V)$
is finite-dimensional. When the finite group $G$ is abelian, this
amounts to the study of Nichols algebras of diagonal type. In this
context, substantial results were reached in \cite{H4}, see also
\cite{AS2, H1, H2, H3}.

For a general finite group $G$, irreducible Yetter-Drinfeld
modules (up to isomorphisms) over $\ku G$ are parametrized by
pairs $(\C, \rho)$ where $\C$ is a conjugacy class of $G$ and
$\rho$ is an irreducible representation of the centralizer $G^s$
of a fixed $s\in \C$. Say $M(\C, \rho)$ is the irreducible
Yetter-Drinfeld module corresponding to $(\C, \rho)$ and
$\toba(\C, \rho)$ is its Nichols algebra. As in \cite{AZ}, we use
the following strategy.

%\smallbreak
{\bf  Strategy.} \emph{ Given $(\C, \rho)$, find a braided
subspace $W$ of $M(\C, \rho)$ of diagonal type. Check if the
dimension of the Nichols algebra $\toba(W)$ is infinite using the
above mentioned results. If so, then necessarily $\dim(\toba(\C,
\rho))=\infty$.}

%\smallbreak
Note that this strategy only involves the determination of a
braided subspace $W$ of $M(\C, \rho)$ of diagonal type. This is an
elementary problem of group theory and its solution does not
requires any knowledge of pointed Hopf algebras or Nichols
algebras.

In this paper, we deal with the groups $G=\sn$. The orbit of $\pi
\in \sn$ is determined by the lengths of the disjoint cycles in
the decomposition of $\pi$. We call $\pi$ \emph{unmixed} if all
those lengths are equal, and \emph{mixed} otherwise.

The main results in this paper are summarized in the following
statement, see below for the unexplained notations.

\begin{maintheorem}\label{mainresults}
Let $\pi \in \mathbb S_{kn} $ be unmixed, say of type $(k^n)$, and $\rho \in \widehat{\mathbb S_{kn}^{\trasp}}$.\\
(A) If $k$ is odd, then $\dim \toba(\C, \rho)= \infty$.\\
\noindent(B) Assume that $k=2$.
\begin{itemize}
\item[(i)] If $n$ is even, then $\dim \toba(\C, \rho)= \infty$.
\item[(ii)] Assume that $n$ is odd, $n > 1$. If $\rho = \chi_{(n)}
\otimes \epsilon$ or $\chi_{(n)} \otimes \sgn$, then the braiding
is negative. Otherwise, $\dim \toba(\C, \rho)= \infty$.
\end{itemize}
\noindent(C) Assume that $k=2r$, with $r > 1$.
\begin{itemize}
\item[(i)] If $n=1$, then $\dim \toba(\C, \rho )= \infty$ if $\rho=\chi_{\omega}$, with $\omega \neq -1$,
and the braiding is negative if $\rho=\chi_{\omega}$, with $\omega
= -1$.
\item[(ii)] Assume that $n>1$. If $\deg \rho>1$, or if $\deg \rho=1$
and $\rho(\pi)\neq -1$, then $\dim \toba(\C, \rho)= \infty$.
Assume that $\deg \rho=1$ and $\rho(\pi)= -1$. If
$\rho=\chi_{r,\dots,r} \otimes \mu$, with $r$ even or odd, or if
$\rho=\chi_{c,\dots,c} \otimes \mu$, with $r$ even and
$c=\frac{r}{2}$ or $\frac{3r}{2}$, then the braiding is negative,
where $\mu=\epsilon$ or $\sgn$; otherwise, $\dim \toba(\C, \rho)=
\infty$.
\end{itemize}
\end{maintheorem}

Part (A) follows from \cite{AZ}, see Lemma \ref{odd}; \cite[Theor.
2.7]{AZ}, Propositions \ref{claim2}, \ref{prop2} and Theorem
\ref{dosn} for part (B); and Proposition \ref{ciclo}, Theorem
\ref{teordegrho>1}, Proposition \ref{propfinal1} and Theorem
\ref{teobra} for part (C). Briefly, for the braided vector spaces
considered in this paper, either the Nichols algebra is
infinite-dimensional or the braiding is negative. The notation in
the Theorem can be found in the pages \pageref{braineg},
\pageref{chi(n)}, \pageref{chi(omega)} and \pageref{ecuacion2r}
for negative braiding, $\chi_{(n)}$, $\chi_{\omega}$ and
$\chi_{r,\dots,r}$, respectively.

We shall consider the mixed case in a subsequent publication;
partial results follow from the unmixed case and \cite[Prop.
2.6]{AZ}.

\section{Preliminaries}\label{conventions}

Our main reference for the classification problem of pointed Hopf
algebras is \cite{AS-cambr}. We denote by $\widehat{G}$ the set of
isomorphism classes of irreducible representations of a finite
group $G$. Consequently, the group of characters of a finite
abelian group $\Gamma$ is denoted $\widehat{\Gamma}$. We shall
often denote a representant of a class in $\widehat{G}$ with the
same symbol as the class itself. If $\rho \in \widehat{G}$, $\deg
\rho$ is the dimension of the vector space $V$ affording $\rho$.
If $V$ is a $\Gamma$-module then $V^\chi$ is the isotypic
component of type $\chi\in \widehat{\Gamma}$. We shall use the
rack notation $g\trid h = ghg^{-1}$, $g, h\in G$. A left comodule
over the group algebra $\ku G$ is the same as a $G$-graded vector
space: if $M$ is a $\ku G$-comodule, then $M = \oplus_{h\in G}
M_h$ where $M_h = \{m\in M: \delta(m) = h\otimes m\}$.

We denote by $\G_n$ the group of $n$-th roots of 1 in $\ku$.

\subsection{Yetter-Drinfeld modules over a finite group}\label{conventionsyd}

Let $G$ be a finite group. A Yetter-Drinfeld module over $G$ is a
left $G$-module and left $\ku G$-comodule $M$ satisfying the
compatibility condition $\delta(g.m) = ghg^{-1} \otimes g.m$, for
all $m\in M_h$, $g, h\in G$. It is well-known that Yetter-Drinfeld
modules over $G$ are completely reducible, and that irreducible
Yetter-Drinfeld modules over $G$ are parameterized by pairs $(\C,
\rho)$ where $\C$ is a conjugacy class in $G$ and $\rho$ is an
irreducible representation of the isotropy subgroup $G^s$ of a
fixed point $s\in \C$ on a vector space $V$.  We denote the
corresponding Yetter-Drinfeld module by $M(\C, \rho)$; a precise
description is as follows. Let $t_1 = s$, \dots, $t_{M}$ be a
numeration of $\C$ and let $g_i\in G$ such that $g_i \trid s =
t_i$ for all $1\le i \le M$. Then  $ M(\C, \rho) = \oplus_{1\le i
\le M}\, g_i\otimes V$. Let $g_iv := g_i\otimes v \in M(\C,\rho)$,
$1\le i \le M$, $v\in V$. If $v\in V$ and $1\le i \le M$, then the
action of $g\in G$ is given by $g\cdot (g_iv) = g_j(\gamma\cdot
v)$, where $gg_i = g_j\gamma$, for some $1\le j \le M$ and
$\gamma\in G^s$, and the coaction is given by $\delta(g_iv) =
t_i\otimes g_iv$. The Yetter-Drinfeld module $M(\C, \rho)$ is a
braided vector space (see below) with braiding
\begin{equation} \label{yd-braiding}
c(g_iv\otimes g_jw) = t_i\cdot(g_jw)\otimes g_iv = g_h(\gamma\cdot
w)\otimes g_iv,
\end{equation}
for any $1\le i,j\le M$, $v,w\in V$, where $t_ig_j = g_h\gamma$
for unique $h$, $1\le h \le M$ and $\gamma \in G^s$. Since $s\in
Z(G^s)$, the Schur Lemma implies that
\begin{equation}\label{schur} s \text{ acts by a scalar $q_{ss}$
on } V.
\end{equation}

Notice that $M(\C, \rho)$ depends on the element $s$ in $\C$ and
$\rho$ in $\widehat{G^s}$. Let $s,\tilde{s} \in \C$ and let $g \in
G$ such that $g s g^{-1}=\tilde{s}$; thus $g G^s
g^{-1}=G^{\tilde{s}}$. Let $\widetilde{\rho}\in
\widehat{G^{\tilde{s}}}$ the pullback of $\rho \in \widehat{G^s}$
via conjugate by $g$. Then $M(\C, \rho)=M(\C,\widetilde{\rho})$;
in particular
\begin{align}\label{image}
\text{the images of $\rho$ and $\widetilde{\rho}$ in $GL(V)$
coincide}.
\end{align}

\subsection{On Nichols algebras}\label{conventionsnichols}

Let $(V, c)$ be a braided vector space, i. e. $V$ is a vector
space and $c: V\otimes V\to V\otimes V$ is an automorphism
satisfying the braid equation $(c\otimes \id)(\id \otimes
c)(c\otimes \id) = (\id \otimes c)(c\otimes \id)(\id \otimes c)$.
Then $\toba(V)$ denotes the Nichols algebra of $V$, see
\cite{AS-cambr}.  The Nichols algebra of a Yetter-Drinfeld module
$M(\C, \rho)$ will be denoted just by $\toba(\C, \rho)$.

\begin{lema}\label{trivialbraiding}
If $W$ is a subspace of $V$ such that $c(W\otimes W) = W\otimes W$
and $\dim \toba(W) =\infty$ then $\dim \toba(V) =\infty$. \qed
\end{lema}

Indeed, $\toba(W) \subset \toba(V)$. A first application of the
Lemma is when there exists $v\in V - 0$ such that $c(v\otimes v )
= v \otimes v$; then $\dim \toba(V) =\infty$. In particular, if $V
= M(\C, \rho)$ and $q_{ss} = 1$, see \eqref{schur}, then
$\dim\toba(\C, \rho) = \infty$.

\medbreak A braided vector space $(V,c)$ is of \emph{diagonal
type} if there exists a basis $v_1, \dots, v_{\theta}$ of $V$ and
non-zero scalars $q_{ij}$, $1\le i,j\le \theta$, such that
$c(v_i\otimes v_j) = q_{ij} v_j\otimes v_i$, for all $1\le i,j\le
\theta$. A braided vector space  $(V,c)$ of diagonal type is of
\emph{Cartan type} if  $q_{ii}\neq 1$ is a root of 1 for all $i$:
$1\le i \le \theta$, and there exists $a_{ij} \in \Z$, $-\ord
q_{ii} < a_{ij} \leq 0$ such that $q_{ij}q_{ji} = q_{ii}^{a_{ij}}$
for all $1\le i\neq j\le \theta$. Set $a_{ii}=2$ for al $1\le i\le
\theta$. Then $(a_{ij})_{1\le i,j\le \theta}$ is a generalized
Cartan matrix.
\begin{theorem}\label{cartantype} (\cite[Th. 4]{H3}, see also
\cite[Th. 1]{AS2}). Let $(V,c)$ be a braided vector space of
Cartan type. Then $\dim \toba(V) < \infty$ if and only if the
Cartan matrix is of finite type. \qed\end{theorem}

\subsection{Abelian subspaces of a braided vector space}\label{cartan-subsection}
Our aim now is to describe a recipe for finding braided subspaces
$W$ of diagonal type of a braided vector space $M(\C, \rho)$.

 Let  $(X, \trid)$ be a rack.  Let $q:
X\times X\to \ku^{\times}$ be a rack 2-cocycle and let $(\ku X,
c_q)$ be the associated braided vector space, that is $\ku X$ is a
vector space with a basis $e_x$, $x\in X$, and $c_q(e_x\otimes
e_y) = q_{x,y} e_{x\trid y}\otimes e_x$, for all $x, y\in X$. Let
us say that a subrack $T$ of $X$  is \emph{abelian} if $i\trid j =
j$ for all $i,j\in T$. If $T$ is an abelian subrack of $X$ then
$\ku T$ is a braided vector subspace of $(\ku X, c_q)$.

\begin{definition}\label{weaklyfinite}
We say that $(\ku X, c_q)$ is \emph{weakly finite} if
$\dim\toba(\ku T) < \infty$ for any abelian subrack $T$ of $X$.
\end{definition}

Our interest is to check when $(\ku X, c_q)$ is not weakly finite,
for then $\dim\toba(\ku X) = \infty$.

We shall say that $(\ku X, c_q)$ is \emph{negative}\label{braineg}
if $q_{ii} = -1$ and $q_{ij}q_{ji} = 1$ for all $i,j\in T$ (hence
$\toba(\ku T)$ is an exterior algebra and $\dim \toba(\ku T) =
2^{\card T}$) and for any abelian subrack $T$ of $X$. This is a
very particular case, but we note that almost all known braided
vector spaces with finite dimensional Nichols algebra that ``do
not come from the abelian case" are negative. See \cite{G2}.

\smallbreak Let now $G$ be a finite group, $\C$ a conjugacy class
in $G$, $\rho\in \widehat{G^s}$ with $s\in \C$ fixed. As in
subsection \ref{conventionsyd}, we fix a numeration $t_1 = s$,
\dots, $t_{M}$ of $\C$ and $g_i\in G$ such that $g_i \trid s =
t_i$ for all $1\le i \le M$. Let $T = \{t_i: i\in I\}$ be an
abelian subrack of $\C$, $I\subset \{1, \dots, M\}$. Let $i, j \in
I$. Then the following are equivalent:
\begin{enumerate}
    \item[(i)] $t_i\trid t_j = t_j$, that is, $t_i$ and $t_j$ commute.
    \item[(ii)]  $\gamma_{ij} := g_j^{-1}t_ig_j\in G^s$.
\end{enumerate}

Let $V$ be the vector space affording $\rho$. For every $v,w\in
V$, we have that
\begin{equation}\label{braiding:abeliansubrack}
c(g_iv\otimes g_jw) =  g_j(\gamma_{ij}\cdot w)\otimes
g_iv\end{equation}
Let $v_1, \dots, v_R$ be simultaneous
eigenvectors of $\gamma_{ij}$, $i, j\in I$. Then
$$
W = \ku-\text{span of } g_i v_j, \quad i\in I, \, 1\le j \le R,
$$
is a braided subspace of diagonal type of dimension $\card T R$.
Note that $R$ depends not only on $T$ but also on the
representation $\rho$; for instance if $\rho$ is a character then
$R = 1 = \dim V$, and $M(\C, \rho)$ is of rack type.

Notice that the action of $G$ on $\C$ (by conjugation) preserves
abelian racks. It is then natural to ask: Are two \emph{maximal}
abelian subracks of $\C$ conjugated by some $g\in G$?

%%%%%%%%%%%%%%%%%%%%%%%%%%%%%%%%%%%%%%%%%%%%%%%%%%%%%%%%%%%%%%%%%%%%%
%%%%%%%%%%%%%%%%%%%%%%%%%%%%%%%%%%%%%%%%%%%%%%%%%%%%%%%%%%%%%%%%%%%%%
%%%%%%%%%%%%%%%%%%%%%%%%%%%%%%%%%%%%%%%%%%%%%%%%%%%%%%%%%%%%%%%%%%%%%

\section{On Nichols algebras over $\sn$}\label{nichols-sn}

\subsection{Notation on the groups $\sn$}\label{conventionss3}

Assume that in the decomposition of $\trasp\in \sn$ as product of
disjoint cycles, there are $m_j$ cycles of length $j$, $1\le j \le
n$. Then the type of  $\trasp \in \sn$ is the symbol $(1^{m_1},
2^{m_2}, \dots, n^{m_n})$; we may omit $j^{m_j}$ when $m_j = 0$.
The conjugacy class $\Oc_\trasp$ of $\trasp$ is the set of all
permutations in $\sn$ with the same type as $\trasp$; we may use
the type as a subscript of a conjugacy class as well. We say that
$\pi$ is \emph{unmixed} if the type of $\pi$ is $(k^n)$, i.e. if
$\pi$ is a product of disjoint cycles of the same length. Let us
assume that $\pi$ is unmixed. It is known that the isotropy
subgroup of $\pi$ satisfies $\mathbb S_{kn}^{\pi}\simeq \Gamma
\rtimes \mathbb S_{n}$, where $\Gamma \simeq (\Z/ k)^{n}$ is
generated by the $k$-cycles in $\pi$ and $\mathbb S_{n}$ permutes
these cycles. This leads us to the representation theory of groups
of the form $\Gamma \rtimes A$, with $\Gamma$ abelian. The
irreducible representations of $\Gamma \rtimes A$ are described as
follows. Let $\chi\in \widehat \Gamma$ and let $A^\chi$ be the
isotropy group with respect to the natural action of $A$ on
$\widehat \Gamma$. Let $\mu\in \widehat{A^\chi}$ and let $\rho$ be
the induced representation $ \rho = \Ind_{\Gamma\rtimes
A^\chi}^{\Gamma\rtimes A} (\chi\otimes \mu)$. Then $\rho$ is
irreducible and any irreducible representation of $\Gamma\rtimes
A$ is isomorphic to one of this form, for unique $\chi$ up to the
action of $A$ and $\mu\in \widehat{A^\chi}$, see \cite[8.2]{S}.

Let $\epsilon$ and $\sgn$ denote the trivial and sign
representation of $\mathbb S_n$, respectively.

\subsection{Nichols algebras corresponding to even
cycles}\label{evencycles}

The next application of Theorem \ref{cartantype} is \cite[Lemma
2.3]{AZ}.

\begin{lema}\label{odd}  If $\dim\toba(\Oc_{\pi}, \rho)< \infty$
then $q_{\pi\pi} = -1$ and $\pi$ has even order.\qed
\end{lema}

In this subsection we consider the case when the type of $\pi$ in
$\mathbb S_{k}$ is $(k)$, with $2 < k$ even. Thus, $\Oc_{\pi}$ is
the set of $k$-cycles in $\mathbb S_{k}$. Fix $\pi=(12\dots k)$;
the isotropy subgroup is $\mathbb S_{k}^{\pi} = \langle \pi
\rangle \simeq \Z_{k}$. If $(j,k)$ denotes the highest common
divisor of $j$ and $k$, then the maximal abelian subrack of
$\Oc_{\pi}$ is
$$
T = \{\pi^j: \, (j,k) = 1\}.
$$
Clearly, $\card T \geq 2$. Let $\omega \in \G_k$, let
$\chi_{\omega}$\label{chi(omega)} be the character of $\mathbb
S_{k}^{\pi}$ defined by $\chi_{\omega}(\pi) = \omega$, let
$M(\Oc_{\pi}, \chi_{\omega})$ be the corresponding Yetter-Drinfeld
module and let $\toba(\Oc_{\pi}, \chi_{\omega})$ be its Nichols
algebra. We conclude from Lemma \ref{odd}:

\begin{prop}\label{ciclo} Let $\trasp\in \mathbb S_k$ of type
$(k)$,  $k$  even. Let  $\omega \in \G_k$. Then $\dim
\toba(\Oc_{\pi}, \chi_{\omega})= \infty$ if $\omega \neq -1$, and
the braiding is negative if $\omega = -1$. \qed
\end{prop}

\subsection{Nichols algebras of orbits with $n$
transpositions}\label{NAntransp}

In this subsection we consider the case when the type of $\pi$ in
$\mathbb S_{2n}$ is $(2^n)$, $n>1$. Thus, $\Oc_{\pi}$ is the set
of permutations in $\mathbb S_{2n}$ that are the product of $n$
disjoint transpositions. Fix $\pi=A_1 \dots A_n$ in $\mathbb
S_{2n}$, with $A_i=(2i-1 \, \, 2i)$. The isotropy subgroup is
$$
\mathbb S_{2n}^{\pi}=\langle A_1,\dots,A_n \rangle \rtimes \langle
B_1,\dots,B_{n-1}\rangle \simeq \Z_2^{n}\rtimes \mathbb S_n,
$$
where $B_j=\left( 2j-1 \quad 2j+1 \right) \left(2j \quad
2j+2\right)$, $1\leq j \leq n-1$; the relations are
\begin{align*}
&A_i^2=\id=B_j^2,  \,\, &&A_iA_j=A_jA_i,  \,\,   &&A_iB_j=B_jA_i \, \text{, $i\neq j, j+1$,}\\
&A_jB_j=B_jA_{j+1}, \,\,  &&B_iB_j=B_jB_i  \, \text{, $|i-j|>1$,}
\,\, &&B_jB_{j\pm 1}B_j=B_{j\pm 1}B_jB_{j\pm 1},
\end{align*}

\bigbreak

\subsubsection{Irreducible representations of $\Z_2^{n}\rtimes
\mathbb S_n$}\label{irrep}

We first list the irreducible representations of $\Z_2^{n}\rtimes
\mathbb S_n$. Let $e_i\in \Z_2^{n}$ be the element with 1 in the
$i$-th component and 0 elsewhere; let $\chi_i \in
\widehat{\Z_2^{n}}$ be given by $\chi_i(e_j)=(-1)^{
\delta_{i,j}}$. The irreducible representations of $\Z_2^{n}$ are
the linear characters $$\chi_{i_1,...,i_k}\label{chi(n)}:=
\chi_{i_1}\dots\chi_{i_k}, \quad 0 \leq k \leq n,\quad 1 \leq i_1
< \cdots <i_k \leq n,$$ where $k=0$ corresponds to the trivial
representation $\chi_{(0)}$ of $\Z_2^{n}$. Let
$\chi_{(k)}:=\chi_{1,\dots,k}$. The action of $\mathbb S_n$ on
$\Z_2^{n}$ induces a natural action of $\mathbb S_n$ on
$\widehat{\Z_2^{n}}$; the orbit and the isotropy subgroup of
$\chi=\chi_{i_1,...,i_k} \in \widehat{\Z_2^{n}}$ are
\begin{align*} \Oc_{\chi}=\{
\chi_{j_1,\dots,j_k} : 1 \leq j_1 <\dots<j_k \leq n \}, \qquad
 \mathbb
S_n^{\chi} \simeq \mathbb S_{n-k} \times \mathbb S_k.
\end{align*}
Thus the characters $\chi_{(k)}$, $0 \leq k \leq n$, are a
complete set of representatives of the orbits in
$\widehat{\Z_2^{n}}$. As discussed in subsection
\ref{conventionss3}, all the irreducible representations of
$\Z_2^{n}\rtimes \mathbb S_n$ are of the form
$$
\rho =\rho_{\chi_{(k)}, \mu}= \Ind_{\Z_2^{n} \rtimes \mathbb
S_n^{\chi_{(k)}}}^{\Z_2^{n} \rtimes \mathbb S_n}
(\chi_{(k)}\otimes \mu), \qquad 0 \leq k \leq n, \quad \mu\in
\widehat{\mathbb S_n^{\chi_{(k)}}}.
$$
There are $\sum_{k=0}^{n} \mathcal P(n-k)\mathcal P(k)$
irreducible representations of $\Z_2^{n}\rtimes \mathbb S_n$,
where $\mathcal P$ is the partition function, but we do not need
to consider all of them.

\begin{obs}\label{claim2bis} If $k$ is even, then
$\rho_{\chi_{(k)}, \mu}(\pi)$ acts by $q_{\pi \pi}=1$, for any
$\mu$ $\in \widehat{\sn^{\chi_{(k)}}}$. Thus
$\dim\toba(\Oc_\trasp, \rho) = \infty$, by Lemma \ref{odd}. So, we
can assume $k$ odd.
\end{obs}

Let $t_1=\pi,\dots,t_M$ be a numeration of $\Oc_{\pi}$, as in
subsection \ref{conventionsyd}; we can assume that the elements
$g_1, \dots, g_M$ satisfying $g_i \trid \pi= t_i$, are
involutions.

\begin{prop}\label{claim2} If $n$ is odd, then the braided vector space associated to
$\chi_{(n)} \otimes \epsilon$ or to $\chi_{(n)} \otimes \sgn$ is
negative. \end{prop}

\pf Assume that $t_i\neq t_j$ commute. We must show that
$q_{ii}=-1$, $q_{jj}=-1$ and $q_{ij}q_{ji}=1$. The first two
conditions are fulfilled because $t_l g_l=g_l \pi$, $1\leq l \leq
M$. For the third, note that $\gamma_{ij} := g_j^{-1}t_ig_j$ and
$\gamma_{ji} := g_i^{-1}t_jg_i$ are in $\snn^{\pi}$, see
subsection \ref{cartan-subsection}; so we can write
\begin{align*}
\gamma_{ij}=A_1^{d_1}\cdots A_n^{d_n}B_{h_1}\cdots B_{h_P}, \qquad
\gamma_{ji}=A_1^{e_1}\cdots A_n^{e_n}B_{l_1}\cdots B_{l_Q},
\end{align*}
where $d_1,\dots , d_n,e_1,\dots ,e_n \in \{0,1\}$. Since
$\gamma_{ij}$, $\gamma_{ji}\in \Oc_{\pi}$, signs of the
permutations $\gamma_{ij}$ and $\gamma_{ji}$ are equal to the sign
of $\pi$, which is $-1$, because $n$ is odd and the sign of every
permutation $B_l$ is $1$, $1\leq l \leq n-1$. This implies that
$d_1 + \cdots + d_n$ and $e_1 + \cdots + e_n$ are odd. Now, since
$t_i g_j= g_j \gamma_{ij}$ and $t_j g_i= g_i \gamma_{ji}$ then
$q_{ij}q_{ji}=\rho(\gamma_{ij}\gamma_{ji})$. We consider the two
cases.

(a) \emph{Assume that $\rho=\chi_{(n)} \otimes \epsilon$}. In this
case, the result follows because
$$\rho(\gamma_{ij}\gamma_{ji})=(-1)^{d_1 + \cdots + d_n}(-1)^{e_1
+ \cdots + e_n}=1.$$

(b) \emph{Assume that $\rho=\chi_{(n)} \otimes \sgn$}. If
$t_i=\pi$, then $\gamma_{ij}=\gamma_{ji}=t_j$, because $g_j$ is an
involution, and the result follows.

We will see that the general case follows from the case $t_i=\pi$.
By definition, $M(\C, \rho)$ is a $\snn$-comodule, with coaction
given by $\delta(g_lv) = t_l\otimes g_lv$, where $V=\ku$ - span of
$v$. Then $ M(\C, \rho) = \oplus_{\tau \in \snn}\, M(\C,
\rho)_{\tau}$, where
$$
M(\C, \rho)_{\tau}:=\{m \in M(\C, \rho) \, : \, \delta(m)=\tau
\otimes m \}.
$$
Clearly, $M(\C, \rho)_{\tau}=g_i V$, if $\tau=t_i$, and $M(\C,
\rho)_{\tau}=0$, if $\tau \notin \Oc_{\pi}$.

Let us call $\widetilde{t_1}:=t_i$; then $M(\C,\rho)\simeq
\Ind_{\mathbb S_{2n}^{\widetilde{t_1}}}^{\mathbb S_{2n}}
\widetilde{V}$, where $\widetilde{V}$ is an irreducible
representation of dimension $1$ of $\mathbb
S_{2n}^{\widetilde{t_1}} \simeq \mathbb S_{2n}^{\pi}$, see
\eqref{image}. Let $\widetilde{t_2}:=t_j$, $\widetilde{g_1}:=\id$
and let $\widetilde{g_2}$ be such that
$\widetilde{g_2}\widetilde{t_1}\widetilde{g_2}^{-1}=\widetilde{t_2}$.
Thus, there exists $\widetilde{v}$ in $\widetilde{V}$ which
satisfies
\begin{align*}
g_i \ku v= M(\C, \rho)_{t_i}= \widetilde{g_1} \ku \widetilde{v},
\qquad \qquad g_j \ku v= M(\C, \rho)_{t_j}= \widetilde{g_2} \ku
\widetilde{v}.
\end{align*}
Let us say $\widetilde{g_1}  \widetilde{v}=\lambda_1 g_i v$ and
$\widetilde{g_2}  \widetilde{v}=\lambda_2 g_j v$. Then
\begin{align*}
c(\widetilde{g_1}  \widetilde{v} \otimes \widetilde{g_2}
\widetilde{v})&=c(\lambda_1 g_i v \otimes \lambda_2 g_j v )=
\lambda_1\lambda_2 c( g_i v \otimes g_j v )\\&=\lambda_1\lambda_2
q_{ij} g_j v \otimes g_i v =q_{ij}\widetilde{g_2}  \widetilde{v}
\otimes \widetilde{g_1} \widetilde{v},
\end{align*}
and on the other hand, $c(\widetilde{g_1}  \widetilde{v} \otimes
\widetilde{g_2} \widetilde{v})=\widetilde{q_{12}} \widetilde{g_l}
\widetilde{v} \otimes \widetilde{g_1} \widetilde{v}$; so,
$q_{ij}=\widetilde{q_{12}}$. Analogously,
$q_{ji}=\widetilde{q_{21}}$. Hence,
$q_{ij}q_{ji}=\widetilde{q_{12}}\widetilde{q_{21}}=1$, and the
result follows. \epf

%%%%%%%%%%%%%%%%%%%%%%%%%%%%%%%%%%%%%%%%%%%%%%%%%%%%%%%%%%%%%%%
%%%%%%%%%%%%%%%%%%%%%%%%%%%%%%%%%%%%%%%%%%%%%%%%%%%%%%%%%%%%%%%
%%%"...Let $\pi=A_1 \dots  A_n$ be in $\mathbb S_{2n}$..."
%%%  LO PASAMOS AL CASO N>4
%%%%%%%%%%%%%%%%%%%%%%%%%%%%%%%%%%%%%%%%%%%%%%%%%%%%%%%%%%%%%%%

\smallbreak We proceed now to consider the different cases
according to the parity of $n$. The case $n=2$ is contained in
\cite{AZ}. We first consider $n=3$, $4$ and then the general cases
$n$ even and $n$ odd.

\subsubsection{Case $n=3$}\label{cason3}

Let $\pi=(1 \, 2)(3 \, 4)(5 \, 6)$ in $\mathbb S_6$. Then
$\Oc_{\pi}$ has $15$ elements and the isotropy subgroup of $\pi$
is
\begin{align*}
\mathbb S_6^{\pi} &= \langle A_1 \!=\! (1 \, 2),A_2 \! = \! (3 \,
4),A_3 \! = \! (5 \, 6)\rangle \rtimes \langle B \! = \! (1 \,
3)(2 \, 4),C \! = \! (1 \, 3 \, 5)(2 \,
4 \, 6)\rangle \\
&\simeq \Z_2^{3}\rtimes \mathbb S_3.
\end{align*}
The defining relations for the generators $A_1,A_2,A_3,B$ and $C$
are $B^2=C^3=1=A_i^2$, $A_i A_j=A_j A_i$, $BCB=C^{-1}$ and
\begin{align*}
BA_1B&=A_2, \qquad & BA_2B&=A_1, \qquad & BA_3B&=A_3,\\
CA_1C^{-1}&=A_2, \qquad & CA_2C^{-1}&=A_3,\qquad &
CA_3C^{-1}&=A_1.
\end{align*}

By section \ref{irrep}, the irreducible representations of
$\mathbb S_6^{\pi}$ are: \\ (1) Four characters $\chi_{\pm,\pm}$
given by $\chi_{\pm,\pm}(A_i)=\pm1$ (the first subindex),
$\chi_{\pm,\pm}(B)=\pm1$ (the second subindex), $\chi_{\pm,\pm}(C)=1$.\\
(2) Two 2-dimensional representations $\theta_{\pm}$ given by
$$\theta_{\pm}(A_i)=\begin{pmatrix}  \pm 1 & 0\\ 0 & \pm 1
\end{pmatrix},\quad  \theta_{\pm}(B)=\begin{pmatrix}  0 & 1\\ 1 & 0
\end{pmatrix},\quad  \theta_{\pm}(C)=\begin{pmatrix}  \omega & 0\\ 0 &
\omega^{-1}
\end{pmatrix},$$
where $\omega \in \ku_3^{\times}$ is a primitive $3$-th root of
1.\\
(3) Four 3-dimensional representations $\phi_{\pm}$, $\psi_{\pm}$
given by
$$\phi_{\pm}(A_1)=\begin{pmatrix}  1 & 0 & 0\\ 0 & -1 &
0\\0 & 0 & 1 \end{pmatrix}, \, \, \phi_{\pm}(A_2)=\begin{pmatrix}
1 & 0 & 0\\ 0 & 1 & 0\\0 & 0 & -1 \end{pmatrix},$$
$$\phi_{\pm}(A_3)=\begin{pmatrix}  -1 & 0 & 0\\
0 & 1 & 0\\0 & 0 & 1\end{pmatrix}, \phi_{\pm}(B)=\begin{pmatrix}
\pm 1 & 0 & 0\\ 0 & 0 & \pm 1 \\0 & \pm 1 & 0
 \end{pmatrix},
\phi_{\pm}(C)=\begin{pmatrix} 0 & 0 & 1\\
1 & 0 & 0\\0 & 1 & 0 \end{pmatrix},$$ and
$$\psi_{\pm}(A_1)=\begin{pmatrix}  -1 & 0 & 0\\ 0 & 1 &
0\\0 & 0 & -1 \end{pmatrix},\psi_{\pm}(A_2)=\begin{pmatrix} -1 & 0
& 0\\ 0 & -1 & 0\\0 & 0 & 1 \end{pmatrix},$$
$$\psi_{\pm}(A_3)\! =\! \begin{pmatrix}  1 & 0 & 0\\
0 & -1 & 0\\0 & 0 &
-1\end{pmatrix},\psi_{\pm}(B) \! =\! \begin{pmatrix} \pm 1 & 0 & 0\\
0 & 0 & \pm 1 \\0 & \pm 1 & 0
 \end{pmatrix},
\psi_{\pm}(C)\! =\! \begin{pmatrix} 0 & 0 & 1\\
1 & 0 & 0\\0 & 1 & 0 \end{pmatrix}.$$

The representations $\chi_{-,\pm}$ are the $\chi_{(n)} \otimes
\epsilon$ and $\chi_{(n)} \otimes \sgn$ in Proposition
\ref{claim2}, thus we can not decide the dimension of their
Nichols algebras. For the others, we have:

\begin{prop}\label{dosdosdos} Let $\trasp\in \mathbb S_6$ with type
$(2^3)$. Let $\rho$ be in $\widehat{\mathbb S_6^{\trasp}}$. If
$\rho\neq\chi_{-,\pm}$ then $\dim\toba(\Oc_\trasp, \rho) =
\infty$.
\end{prop}

\pf We can assume that $\pi=A_1A_2A_3=(12)(34)(56)$. Let $\rho \in
\widehat{\sn^{\pi}}$ and $V$ the vector space affording $\rho$. We
look for a braided subspace of diagonal type of $M(\Oc_{2,2,2},
\rho )$. Set $\alpha=(12)(35)(46)$, $\beta=(12)(36)(45)$ in
$\Oc_{2,2,2}$; if $\sigma_1=\id$, $\sigma_2=(4 \, 5)$,
$\sigma_3=(4 \, 6)$ then
$$\sigma_1 \trid \pi= \pi, \quad \sigma_2 \trid \pi= \alpha, \quad \sigma_3 \trid \pi= \beta.$$

Let $\sigma_jv := \sigma_j \otimes v$, $v\in V$, $1\leq j \leq 3$.
The coaction is given by $\delta(\sigma_j v) = \sigma_j \trid \pi
\otimes \sigma_j v$; we need the action of the elements $\pi$,
$\alpha$, $\beta$, which is
\begin{align*}
\pi\cdot \sigma_1 v &= \sigma_1 \rho(\pi)(v), &\quad \pi\cdot
\sigma_2
v &= \sigma_2 \rho(\alpha)(v), &\quad \pi\cdot \sigma_3 v &= \sigma_3 \rho(\beta)(v), & \\
\alpha\cdot \sigma_1 v &= \sigma_1 \rho(\alpha)(v), &\quad
\alpha\cdot \sigma_2
v &= \sigma_2 \rho(\pi)(v), &\quad \alpha\cdot \sigma_3 v &= \sigma_3 \rho(\alpha)(v), & \\
\beta\cdot \sigma_1 v &= \sigma_1 \rho(\beta)(v), &\quad
\beta\cdot \sigma_2 v &= \sigma_2 \rho(\beta)(v), &\quad
\beta\cdot \sigma_3 v &= \sigma_3 \rho(\pi)(v).
\end{align*}

Hence the braiding is given,  by
\begin{align*}
c(\sigma_1 v \otimes \sigma_2 w) &= \sigma_2 \rho(\alpha)(w)
\otimes \sigma_1 v, \qquad c(\sigma_1 v \otimes \sigma_3 w) &=
\sigma_3
\rho(\beta)(w) \otimes \sigma_1 v,\\
c(\sigma_2 v \otimes \sigma_1 w) &= \sigma_1 \rho(\alpha)(w)
\otimes \sigma_2 v, \qquad c(\sigma_2 v \otimes \sigma_3 w) &=
\sigma_3 \rho(\alpha)(w) \otimes \sigma_2 v, \\ c(\sigma_3 v
\otimes \sigma_1 w) &= \sigma_1 \rho(\beta)(w) \otimes \sigma_3 v,
\qquad c(\sigma_3 v \otimes \sigma_2 w) &= \sigma_2 \rho(\beta)(w)
\otimes \sigma_3 v,
\end{align*}
and $c(\sigma_j v \otimes \sigma_j w) = (\sigma_j\trid \pi)\cdot
\sigma_jw \otimes \sigma_j v = \sigma_j \rho(\pi)(w) \otimes
\sigma_j v = -\sigma_j w \otimes \sigma_j v $, for all $1\leq j
\leq 3$ and $v,w\in V$.

We now consider the different possibilities for $\rho$. If
$\rho=\chi_{+,\pm}$, $\theta_+$ or $\psi_{\pm}$ then $q_{\pi
\pi}=1$ and $\dim \toba(\Oc_{2,2,2}, \rho) = \infty$, by Lemma
\ref{trivialbraiding}.

If $\rho=\theta_-$ then $\rho(\alpha)=\rho(\beta)=
\begin{pmatrix}  0 & -\omega^{-1}\\ -\omega &  0
\end{pmatrix}$. Choose $v_1 = \begin{pmatrix} 1 \\
-\omega \end{pmatrix}$ and $v_2 =
\begin{pmatrix}  1 \\ \omega\end{pmatrix}$. Hence $\rho(\alpha)(v_1) =\rho(\beta)(v_1)= v_1$,
$\rho(\alpha)(v_2) =\rho(\beta)(v_2)= -v_2$. Therefore the
braiding is diagonal of Cartan type in the basis
$$w_1 = \sigma_1v_1,\,\, w_2 = \sigma_1v_2, \,\, w_3 =
\sigma_2v_1,\,\, w_4 = \sigma_2v_2,\,\,  w_5 = \sigma_3v_1,\,\,
w_6 = \sigma_3v_2.$$ The corresponding Dynkin diagram is
$A_5^{(1)}$, which is affine. By Theorem \ref{cartantype}, $\dim
\toba(\Oc_{2,2,2}, \theta_-) = \infty$ .

\medbreak Assume now that $\rho=\phi_+$. Then
$\phi_+(\alpha)=\small{\begin{pmatrix} 0 & 0 & 1\\ 0 & -1 & 0\\1 &
0 & 0
\end{pmatrix}}$,\newline$\phi_+(\beta)=\small{\begin{pmatrix} 0 & 0 & -1\\ 0 &
-1 & 0\\-1 & 0 & 0\end{pmatrix}}$. Choose
$v_1=\small{\begin{pmatrix}0 \\ 1 \\ 0
\end{pmatrix}}$, $v_2=\small{\begin{pmatrix}1 \\ 0 \\ 1
\end{pmatrix}}$ and $v_3=\small{\begin{pmatrix}1 \\ 0 \\ -1
\end{pmatrix}}$. Thus,
\begin{align*}
&\rho(\alpha)(v_1)=-v_1, \qquad &\rho(\alpha)(v_2)&=v_2, \qquad
&\rho(\alpha)(v_3)&=-v_3,\\
& \rho(\beta)(v_1)=-v_1, \qquad & \rho(\beta)(v_2)&=-v_2, \qquad &
\rho(\beta)(v_3)&=v_3.
\end{align*}
Hence the braiding is diagonal of Cartan type in the basis
\begin{align*}
&w_1 = \sigma_1v_1, \quad & w_2&= \sigma_1v_2, \quad & w_3 &=\sigma_1v_3,\\
&w_4 = \sigma_2v_1, \quad & w_5&= \sigma_2v_2, \quad & w_6 &=\sigma_2v_3,\\
&w_7 = \sigma_3v_1, \quad & w_8&= \sigma_3v_2, \quad & w_9 &=
\sigma_3 v_3;
\end{align*}
this implies that the corresponding Dynkin diagram contains the
affine Dynkin diagram $A_2^{(1)}$. By Theorem \ref{cartantype},
$\dim \toba(\Oc_{2,2,2}, \phi_+) = \infty$. Finally, the case
$\rho=\phi_-$ is similar. \epf

\subsubsection{Case $n=4$}\label{cason4} Let $\pi=(1 \, 2)(3 \,
4)(5 \, 6)(7 \, 8)$ in $\mathbb S_8$. The isotropy subgroup of
$\pi$ is $ \mathbb S_8^{\pi}=<A_1,A_2,A_3,A_4> \rtimes
<B_1,B_2,B_3> \simeq \Z_2^{4}\rtimes \mathbb S_4$, where $A_1=(1
\, 2)$, $A_2=(3 \, 4)$, $A_3=(5 \, 6)$, $A_4=(7 \, 8)$, and
$$
B_1=(1 \, 3)(2 \, 4) \quad B_2=(3 \, 5)(4 \, 6) \quad B_3=(5 \,
7)(6 \, 8).
$$
By section \ref{irrep}, there are $20$ irreducible representations
of $\mathbb S_8^{\pi}$, but by Remark \ref{claim2bis} we only need
to consider $6$ of them. They are
\begin{align*}
\rho_1&= \Ind (\chi_{(1)}\otimes \epsilon), \quad & \rho_2&=\Ind
(\chi_{(1)}\otimes \sgn), \quad & \rho_3&=\Ind (\chi_{(1)}\otimes \theta),\\
\rho_4&=\Ind (\chi_{(3)}\otimes \epsilon), \quad & \rho_5&=\Ind
(\chi_{(3)}\otimes \sgn), \quad  & \rho_6&=\Ind (\chi_{(3)}\otimes
\theta),
\end{align*}
where $\Ind$ means $\Ind_{\Z_2^{4} \rtimes \mathbb S_3}^{\Z_2^{4}
\rtimes \mathbb S_4}$ and $\theta$ is the standard representation
of $\mathbb S_3$. With this notation, we can state the following
result.

\begin{prop}\label{dos4} Let $\trasp\in \mathbb S_8$ of type
$(2^4)$. Then $\dim\toba(\Oc_\trasp, \rho) = \infty$ for any $\rho
\in \widehat{\mathbb S_8^{\trasp}}$. \end{prop}

\pf We can assume that $\pi=(1 \, 2)(3 \, 4)(5 \, 6)(7 \, 8)$. We
shall later prove that if $i=1,2,4$ or $5$ then
$\dim\toba(\Oc_\trasp, \rho_i) = \infty$, for any $n\geq 4$, see
Lemmata \ref{lema1es} and \ref{leman-1es}. We check now the
remaining $i=3,6$. It is clear that $\pi$,
$\alpha=(12)(34)(57)(68)$, $\beta=(12)(34)(58)(67)$ are in
$\Oc_{\pi}$ and they satisfy the same relations as in the proof of
Proposition \ref{dosdosdos}, with $\sigma_1=\id$, $\sigma_2=(6 \,
7)$, $\sigma_3=(6 \, 8)$. Let us consider $\rho_3$. In an
appropriate basis, we have

\begin{align*}
\rho_3(\alpha)&=\small{\begin{pmatrix} 0 & 0 & 0 & 0 & -\omega^{-1} & 0 & 0 & 0\\
0 & 0 & 0 & 0 &  0 & -\omega^{-1} &  0 & 0\\
0 & 0 & 0 & 1 &  0 & 0 &  0 & 0\\
0 & 0 & 1 & 0 &  0 & 0 &  0 & 0\\
-\omega & 0 & 0 & 0 &  0 & 0 &  0 & 0\\
0 & -\omega & 0 & 0 &  0 & 0 &  0 & 0\\
0 & 0 & 0 & 0 &  0 & 0 &  0 & 1\\
0 & 0 & 0 & 0 &  0 & 0 &  1 & 0\\
\end{pmatrix}},\end{align*}
\begin{align*}
\rho_3(\beta)&=\small{\begin{pmatrix} 0 & 0 & 0 & 0 & -\omega^{-1} & 0 & 0 & 0\\
0 & 0 & 0 & 0 &  0 & -\omega^{-1} &  0 & 0\\
0 & 0 & 0 & -1 &  0 & 0 &  0 & 0\\
0 & 0 & -1 & 0 &  0 & 0 &  0 & 0\\
-\omega & 0 & 0 & 0 &  0 & 0 &  0 & 0\\
0 & -\omega & 0 & 0 &  0 & 0 &  0 & 0\\
0 & 0 & 0 & 0 &  0 & 0 &  0 & -1\\
0 & 0 & 0 & 0 &  0 & 0 &  -1 & 0\\
\end{pmatrix}},
\end{align*}
where $\omega \in \ku_3^{\times}$ is a primitive $3$-th root of
the unity. It is easy to see that
\begin{align*} &v_1=e_1+ \omega
e_5, \quad & v_2&=e_2+ \omega e_6, \quad v_3=e_3+
e_4, \quad v_4=e_3- e_4,\\
&v_5=e_1- \omega e_5, \quad & v_6&=e_2- \omega e_6, \quad v_7=e_7+
e_8, \quad v_8=e_7- e_8
\end{align*}
are eigenvectors of eigenvalues $1$ or $-1$. In particular
\begin{align*}
\rho_3(\alpha)v_7=v_7, \quad \rho_3(\alpha)v_8=-v_8, \quad
\rho_3(\beta)v_7=-v_7, \quad \rho_3(\beta)v_8=v_8.
\end{align*}
So, in the basis $w_1 = \sigma_1v_7$, $w_2= \sigma_1v_8$, $ w_3
=\sigma_2v_7$, $w_4 = \sigma_2v_8$, $w_5= \sigma_3v_7$, $w_6
=\sigma_3v_8$, the braiding is diagonal of Cartan type. The
corresponding Dynkin diagram is not connected; its connected
components are $\{1, 4, 6\}$ and $\{2, 3, 5\}$, each of them
supporting the affine Dynkin diagram $A_2^{(1)}$. Then $\dim
\toba(\Oc_{\pi}, \rho_3) = \infty$ by Theorem \ref{cartantype}.

Finally, if $\rho=\rho_6$ we have
$$
\rho_6(\alpha)= \small{\begin{pmatrix} 0 & 1 & 0 & 0 & 0 & 0 & 0 & 0\\
1 & 0 & 0 & 0 &  0 & 0 &  0 & 0\\
0 & 0 & 0 & 0 &  0 & 0 &  -\omega^{-1} & 0\\
0 & 0 & 0 & 0 &  0 & 0 &  0 & -\omega^{-1}\\
0 & 0 & 0 & 0 &  0 & 1 &  0 & 0\\
0 & 0 & 0 & 0 &  1 & 0 &  0 & 0\\
0 & 0 & -\omega & 0 &  0 & 0 &  0 & 0\\
0 & 0 & 0 & -\omega &  0 & 0 &  0 & 0\\
\end{pmatrix}},
$$
$$
\rho_6(\beta)=\small{\begin{pmatrix} 0 & -1 & 0 & 0 & 0 & 0 & 0 & 0\\
-1 & 0 & 0 & 0 &  0 &  &  0 & 0\\
0 & 0 & 0 & 0 &  0 & 0 &  -\omega^{-1} & 0\\
0 & 0 & 0 & 0 &  0 & 0 &  0 & -\omega^{-1}\\
0 & 0 & 0 & 0 &  0 & -1 &  0 & 0\\
0 & 0 & 0 & 0 &  -1 & 0 &  0 & 0\\
0 & 0 & -\omega & 0 &  0 & 0 &  0 & 0\\
0 & 0 & 0 & -\omega &  0 & 0 &  0 & 0\\
\end{pmatrix}}.
$$
The following are eigenvectors of eigenvalues $1$ or $-1$
\begin{align*} &v_1=e_1+e_2, \quad & v_2&=e_1-e_2, \quad v_3=e_3+\omega
e_7, \quad v_4=e_4+\omega e_8,\\
&v_5=e_5+ e_6, \quad & v_6&=e_5- e_6, \quad v_7=e_3- \omega e_7,
\quad v_8=e_4- \omega e_8.
\end{align*}
In particular
\begin{align*}
\rho_6(\alpha)v_1=v_1, \quad \rho_6(\alpha)v_2=-v_2, \quad
\rho_6(\beta)v_1=-v_1, \quad \rho_6(\beta)v_2&=v_2.
\end{align*}
So, in the basis $w_1 = \sigma_1v_1$, $w_2= \sigma_1v_2$, $w_3
=\sigma_2v_1$, $w_4 = \sigma_2v_2$, $w_5= \sigma_3v_1$, $w_6
=\sigma_3v_2$, the braiding is diagonal of Cartan type. This is
similar to the case of $\rho_3$ interchanging the roles of $v_7$
by $v_1$ and $v_8$ by $v_2$. Then $\dim \toba(\Oc_{\pi}, \rho_6) =
\infty$ by Theorem \ref{cartantype}. \epf

\subsubsection{Case $n$ general}\label{ngral}

We now begin the analysis of the general case. We prove two
lemmata.

\begin{lema}\label{lema1es}
If $\rho$ is either $\Ind_{\Z_2^{n} \rtimes \mathbb
S_n^{\chi_{(1)}}}^{\Z_2^{n} \rtimes \mathbb S_n}
(\chi_{(1)}\otimes \epsilon)$ or $\Ind_{\Z_2^{n} \rtimes \mathbb
S_n^{\chi_{(1)}}}^{\Z_2^{n} \rtimes \mathbb S_n}
(\chi_{(1)}\otimes \sgn)$ then $\dim\toba(\Oc_\trasp, \rho) =
\infty$, for any $n\geq 2$.
\end{lema}
\pf Let $\sigma_1=\id$, $\sigma_2=(2n-2 \, 2n-1)$ and
$\sigma_3=(2n-2 \, 2n)$. We define $\alpha:=\sigma_2 \trid \pi$
and $\beta= \sigma_3 \trid \pi$, so
\begin{align*}
\alpha&=A_1A_2\dots A_{n-2}(2n-3 \, 2n-1)(2n-2 \, 2n), \\
\beta&= A_1A_2\dots A_{n-2}(2n-3 \, 2n)(2n-2 \, 2n-1);
\end{align*}
we set $T=\{ \pi,\alpha,\beta  \}$. It is straightforward to check
that
\begin{equation}\label{rel1}
\alpha= \pi A_{n-1}A_{n}B_{n-1}=\beta A_{n-1}A_{n}, \qquad
\beta=\pi B_{n-1}.
\end{equation}

Now we proceed as in the proof of Proposition \ref{dosdosdos}. Let
$\sigma_jv := \sigma_j \otimes v$, $v\in V$, $1\leq j \leq 3$. The
coaction is given by $\delta(\sigma_j v) = \sigma_j \trid \pi
\otimes \sigma_j v$; the action of the elements $\pi$, $\alpha$,
$\beta$ is
\begin{align*}
\pi\cdot \sigma_1 v &= \sigma_1 \rho(\pi)(v), &\quad \pi\cdot
\sigma_2
v &= \sigma_2 \rho(\alpha)(v), &\quad \pi\cdot \sigma_3 v &= \sigma_3 \rho(\beta)(v), & \\
\alpha\cdot \sigma_1 v &= \sigma_1 \rho(\alpha)(v), &\quad
\alpha\cdot \sigma_2
v &= \sigma_2 \rho(\pi)(v), &\quad \alpha\cdot \sigma_3 v &= \sigma_3 \rho(\alpha)(v), & \\
\beta\cdot \sigma_1 v &= \sigma_1 \rho(\beta)(v), &\quad
\beta\cdot \sigma_2 v &= \sigma_2 \rho(\beta)(v), &\quad
\beta\cdot \sigma_3 v &= \sigma_3 \rho(\pi)(v).
\end{align*}
Hence the braiding is given by
\begin{align*}
c(\sigma_1 v \otimes \sigma_2 w) &= \sigma_2 \rho(\alpha)(w)
\otimes \sigma_1 v, \qquad c(\sigma_1 v \otimes \sigma_3 w) &=
\sigma_3
\rho(\beta)(w) \otimes \sigma_1 v,\\
c(\sigma_2 v \otimes \sigma_1 w) &= \sigma_1 \rho(\alpha)(w)
\otimes \sigma_2 v, \qquad c(\sigma_2 v \otimes \sigma_3 w) &=
\sigma_3 \rho(\alpha)(w) \otimes \sigma_2 v, \\ c(\sigma_3 v
\otimes \sigma_1 w) &= \sigma_1 \rho(\beta)(w) \otimes \sigma_3 v,
\qquad c(\sigma_3 v \otimes \sigma_2 w) &= \sigma_2 \rho(\beta)(w)
\otimes \sigma_3 v,
\end{align*}
and $c(\sigma_j v \otimes \sigma_j w) = (\sigma_j\trid \pi)\cdot
\sigma_jw \otimes \sigma_j v = \sigma_j \rho(\pi)(w) \otimes
\sigma_j v = -\sigma_j w \otimes \sigma_j v$, for all $1\leq j
\leq 3$ and $v,w\in V$.

Let us consider $\rho=\Ind_{\Z_2^{n} \rtimes \mathbb
S_n^{\chi_{(1)}}}^{\Z_2^{n} \rtimes \mathbb S_n}
(\chi_{(1)}\otimes \epsilon)$. The vector space affording $\rho$
has dimension $n$. It is easy to see that for every $i$, $1\leq i
\leq n$, the matrix $\rho(A_i)$ is diagonal with
$(\rho(A_i))_{i,i}=-1$ and $1$ elsewhere; while
$$
\rho(B_{n-1})=\begin{pmatrix}
\Id_{n-2} &  &  \\
 & 0 & 1 \\
 & 1 & 0 \\
\end{pmatrix}.
$$
Therefore, we have $\rho(\pi)=-\Id$ and, by (\ref{rel1}),
\begin{align*}
\rho(\alpha)=\begin{pmatrix}
-\Id_{n-2} &  & \\
 & 0 & 1 \\
 & 1 & 0 \\
\end{pmatrix}, \quad
\rho(\beta)=\begin{pmatrix}
-\Id_{n-2} &  & \\
 & 0 & -1 \\
 & -1 & 0 \\
\end{pmatrix}.
\end{align*}
Choose $v_i=e_i$, $1\leq i \leq n-2$, $v_{n-1}=e_{n-1}+e_n$ and
$v_n=e_{n-1}-e_n$. Hence
\begin{align*}
\rho(\alpha) v_i&=-v_i, \quad & \rho(\alpha) v_{n-1}&=v_{n-1}, \quad &\rho(\alpha)v_n&=-v_n,\\
\rho(\beta) v_i&=-v_i, \quad &\rho(\beta) v_{n-1}&=-v_{n-1}, \quad
&\rho(\beta)v_n&=v_n,
\end{align*}
with $1 \leq i \leq n-2$. Then the braiding is diagonal of Cartan
type in the basis $\mathcal{B}=\{\sigma_j v_i\}$, $j=1,2,3$,
$1\leq i \leq n$. The corresponding Dynkin diagram is not of
finite type because it contains the affine Dynkin diagram
$A_2^{(1)}$. By Theorem \ref{cartantype}, $\dim \toba(\Oc_{\pi},
\rho) = \infty$.

Finally, if $\rho=\Ind_{\Z_2^{n} \rtimes \mathbb
S_n^{\chi_{(1)}}}^{\Z_2^{n} \rtimes \mathbb S_n}
(\chi_{(1)}\otimes \sgn)$, $\rho(A_i)$ are as before and
$$
\rho(B_{n-1})=\begin{pmatrix}
-\Id_{n-2} &  & \\
 & 0 & 1 \\
& 1 & 0 \\
\end{pmatrix}.
$$
Then,
\begin{align*}
\rho(\alpha)=\begin{pmatrix}
\Id_{n-2} &  & \\
 & 0 & 1 \\
 & 1 & 0 \\
\end{pmatrix}, \quad
\rho(\beta)=\begin{pmatrix}
\Id_{n-2} &  &\\
 & 0 & -1 \\
 & -1 & 0 \\
\end{pmatrix}.
\end{align*}
Let $v_i$ be as before; hence
\begin{align*}
\rho(\alpha) v_i&=v_i, \quad & \rho(\alpha) v_{n-1}&=v_{n-1}, \quad &\rho(\alpha)v_n&=-v_n,\\
\rho(\beta) v_i&=v_i, \quad &\rho(\beta) v_{n-1}&=-v_{n-1}, \quad
&\rho(\beta)v_n&=v_n,
\end{align*}
with $1 \leq i \leq n-2$. Then the braiding is diagonal of Cartan
type in the basis $\mathcal{B}$. The corresponding Dynkin diagram
is not of finite type because it contains the affine diagram
$A_2^{(1)}$. By Theorem \ref{cartantype}, $\dim \toba(\Oc_{\pi},
\rho) = \infty$. \epf

\begin{lema}\label{leman-1es}
If $\rho$ is either $\Ind_{\Z_2^{n} \rtimes \mathbb
S_n^{\chi_{(n-1)}}}^{\Z_2^{n} \rtimes \mathbb S_n}
(\chi_{(n-1)}\otimes \epsilon)$ or \newline $\Ind_{\Z_2^{n}
\rtimes \mathbb S_n^{\chi_{(n-1)}}}^{\Z_2^{n} \rtimes \mathbb S_n}
(\chi_{(n-1)}\otimes \sgn)$ then $\dim\toba(\Oc_\trasp, \rho) =
\infty$, for any $n\geq 2$.
\end{lema}

\pf Let $\sigma_j$, $\alpha$ and $\beta$ be as in the proof of
Lemma \ref{lema1es}.

If $\rho=\Ind_{\Z_2^{n} \rtimes \mathbb
S_n^{\chi_{(n-1)}}}^{\Z_2^{n} \rtimes \mathbb S_n}
(\chi_{(n-1)}\otimes \epsilon)$, then for every $i$, $1\leq i \leq
n$, the matrix $\rho(A_i)$ is diagonal with
$(\rho(A_i))_{n-i+1,n-i+1}=1$ and $-1$ elsewhere; while
$$
\rho(B_{n-1})=\begin{pmatrix}
0 & 1 & \\
1 & 0 &  \\
  &  & \Id_{n-2} \\
\end{pmatrix}.
$$
Therefore, we have $\rho(\pi)=-\Id$ and, by (\ref{rel1}),
\begin{align*}
\rho(\alpha)=\begin{pmatrix}
0 & -1 & \\
-1 & 0 & \\
  &  & \Id_{n-2} \\
\end{pmatrix}, \quad
\rho(\beta)=\begin{pmatrix}
0 & -1 & \\
-1 & 0 & \\
  &  & -\Id_{n-2} \\
\end{pmatrix}.
\end{align*}
Choose $v_1=e_1+e_2$, $v_2=e_1-e_2$ and $v_i=e_i$, $3\leq i \leq
n$. Hence
\begin{align*}
\rho(\alpha) v_1&=-v_1, \quad & \rho(\alpha) v_2&=v_2, \quad &\rho(\alpha)v_i&=v_i,\\
\rho(\beta) v_1&=-v_1, \quad  &\rho(\beta) v_2&=v_2, \quad &
\rho(\beta)v_i&=-v_i,
\end{align*}
with $3 \leq i \leq n$. Then the braiding is diagonal of Cartan
type in the basis $\mathcal{B}=\{\sigma_j v_i\}$, $j=1,2,3$,
$1\leq i \leq n$. The corresponding Dynkin diagram is not of
finite type because it contains the affine Dynkin diagram
$A_5^{(1)}$. By Theorem \ref{cartantype}, $\dim \toba(\Oc_{\pi},
\rho) = \infty$.

\medbreak Finally, if $\rho=\Ind_{\Z_2^{n} \rtimes \mathbb
S_n^{\chi_{(n-1)}}}^{\Z_2^{n} \rtimes \mathbb S_n}
(\chi_{(n-1)}\otimes \sgn)$ then the matrices $\rho(A_i)$ are the
same as in the previous case and
$$ \rho(B_{n-1})=\begin{pmatrix}
0 & 1 & \\
1 & 0 &  \\
  &  & -\Id_{n-2} \\
\end{pmatrix}.
$$
Then,
\begin{align*}
\rho(\alpha)=\begin{pmatrix}
0 & -1 & \\
-1 & 0 & \\
 &  & -\Id_{n-2} \\
\end{pmatrix}, \quad
\rho(\beta)=\begin{pmatrix}
0 & -1 & \\
-1 & 0 & \\
  &  & \Id_{n-2} \\
\end{pmatrix}.
\end{align*}

Let $v_i$ as before; hence
\begin{align*}
\rho(\alpha) v_1&=-v_1, \quad & \rho(\alpha) v_2&=v_2, \quad &\rho(\alpha)v_i&=-v_i,\\
\rho(\beta) v_1&=-v_1, \quad  &\rho(\beta) v_2&=v_2, \quad &
\rho(\beta)v_i&=v_i,
\end{align*}
with $3 \leq i \leq n$. Then the braiding is diagonal of Cartan
type in the basis $\mathcal{B}$. The corresponding Dynkin diagram
is not of finite type because it contains the affine Dynkin
diagram $A_5^{(1)}$. By Theorem \ref{cartantype}, $\dim
\toba(\Oc_{\pi}, \rho) = \infty$. \epf

\smallbreak Notice that if $n$ is odd, then the analog of Lemma
\ref{leman-1es} follows from Remark \ref{claim2bis}. Next we set
some notation separately for the cases $n$ even and $n$ odd.

\bigbreak

\textbf{Notation in case $n$ even}. We suppose $n=2L$ and for
every $l$, with $1\leq l \leq L$, we define $\sigma_l^+:=(4l-2
\,\, 4l-1)$, $\sigma_l^-:=(4l-2 \,\, 4l)$ and
\begin{align*}
\alpha_l:=\sigma_l^+ \trid \pi, \qquad \beta_l:=\sigma_l^- \trid
\pi.
\end{align*}
That is, if $\pi=A_1A_2\dots A_{2l-1}A_{2l} \dots A_{2L-1}A_{2L}$,
then
\begin{align*}
\alpha_l&=A_1A_2\dots (4l-3 \,\, 4l-1)(4l-2 \,\, 4l) \dots A_{2L-1}A_{2L}, \\
\beta_l&=A_1A_2\dots (4l-3 \,\, 4l)(4l-2 \,\, 4l-1) \dots
A_{2L-1}A_{2L}.
\end{align*}
It is easy to see that $\sigma_l^{\pm}
\sigma_h^{\pm}=\sigma_h^{\pm}  \sigma_l^{\pm}$, for all $l,h$
distinct. Let $T$ be the set
$$
T=\{ \sigma_{l_k}^{\pm} \trid (\cdots( \sigma_{l_1}^{\pm} \trid
\pi)\cdots) \, : \,  1\leq k \leq L , 1\leq l_1 <  \cdots < l_k
\leq L\} \cup \{\pi\}.
$$
Note that $\sigma_{l_k}^{\pm} \trid (\cdots( \sigma_{l_1}^{\pm}
\trid \pi)\cdots)=(\sigma_{l_k}^{\pm}\cdots \sigma_{l_1}^{\pm})
\trid \pi$. Let $T\! = \! \{\pi_0 \! = \! \pi, \pi_1,\dots, \pi_N
\}$ be a numeration of $T$; we call $\sigma_0=\id$ and $\sigma_j$
the element $\sigma_l^{\pm}$ such that $\sigma_j \trid \pi=\pi_j$,
$1\leq j \leq N$.

\smallbreak Straightforward computations imply the following.

\begin{lema}\label{lesig}
For all $l$, $1\leq l \leq L$, we have
\begin{itemize}
\item[(i)] $\sigma_l^-  \sigma_l^+ A_{2l-1}A_{2l}\sigma_l^+  \sigma_l^-=\sigma_l^+
A_{2l-1}A_{2l}\sigma_l^+$.
\item[(ii)] $\sigma_l^+  \sigma_l^- A_{2l-1}A_{2l}\sigma_l^-  \sigma_l^+=\sigma_l^-
A_{2l-1}A_{2l}\sigma_l^-$.
\item[(iii)] $\sigma_l^{\pm} A_{2l-1}A_{2l}\sigma_l^{\pm} A_{2l-1}A_{2l}= \sigma_l^{\mp}
A_{2l-1}A_{2l}\sigma_l^{\mp}$.
\item[(iv)] $A_{2l-1}A_{2l} \sigma_l^{\pm} A_{2l-1}A_{2l}\sigma_l^{\pm} = \sigma_l^{\mp}
A_{2l-1}A_{2l}\sigma_l^{\mp}$.
\item[(v)] $\alpha_l\beta_l=A_{2l-1}A_{2l}=\beta_l\alpha_l$.\qed
\end{itemize}
\end{lema}

\bigbreak

\begin{lema}\label{prop1}
(1) If $1\leq l \leq L$, then
\begin{align*}
\alpha_l=\pi A_{2l-1}A_{2l} B_{2l-1}  =\beta_l A_{2l-1}A_{2l},
\qquad \beta_l=\pi B_{2l-1}.
\end{align*}
(2) $T \subseteq \Oc_\trasp \cap \snn^{\pi}$.\\
(3) For every $i,j$, $0\leq i,j \leq N$, there exists $k$,
$0\leq k \leq N$ such that $\pi_i \sigma_j=\sigma_j \pi_k$.\\
(4) $T$ is abelian.
\end{lema}

\pf (1) is obvious. (2) Clearly, $T \subseteq \Oc_\trasp$. To see
that $T \subseteq \snn^{\pi}$ we need to prove that
$\pi_i=(\sigma_{l_k}^{\pm}\cdots \sigma_{l_1}^{\pm}) \trid \pi$ in
$\snn^{\pi}$. This is clear for $k=1$; then it follows by
induction on $k$.

(3) Note that if $i=0$ then $k=j$; if $i=j$ then $k=0.$, etc. Fix
$\pi_i=(\sigma_{l_k}^{\pm}\cdots \sigma_{l_1}^{\pm}) \trid \pi$,
with $1\leq l_1 <  \cdots < l_k \leq L$; suppose
$\sigma_j=\sigma_{h_M}^{\pm}\cdots \sigma_{h_1}^{\pm}$, with
$1\leq h_1 < \cdots < h_M \leq L$. Then
$$
\sigma_j \pi_i \sigma_j= \sigma_{l_k}^{\pm}\cdots
\sigma_{l_1}^{\pm} \sigma_{h_M}^{\pm}\cdots \sigma_{h_1}^{\pm} \pi
\sigma_{h_1}^{\pm}\cdots \sigma_{h_M}^{\pm}
\sigma_{l_1}^{\pm}\cdots \sigma_{l_k}^{\pm}.
$$
If $l_r \neq h_s$ for all $r,s$ then $\pi_k:=\sigma_j \pi_i
\sigma_j$ is in $T$. If $l_r=h_s$ for many $r,s$ we have that in
the expression $\sigma_j \pi_i \sigma_j$ the factors
$\sigma_{l_r}^{\pm}$ and $\sigma_{h_s}^{\pm}$ cancel mutually
while for the factors $\sigma_{l_r}^{\pm}$ and
$\sigma_{h_s}^{\mp}$ we use the Lemma \ref{lesig}(i),(ii) and we
have
$$
\sigma_j \pi_i \sigma_j= \cdots \sigma_{l_r}^{\pm}
\sigma_{h_s}^{\mp} \pi  \sigma_{h_s}^{\mp} \sigma_{l_r}^{\pm}
\cdots = \cdots \sigma_{h_s}^{\mp} \pi \sigma_{h_s}^{\mp} \cdots.
$$
Therefore $\sigma_j \pi_i \sigma_j$ is in $T$.

(4) We need  to prove $\pi_i\pi_j=\pi_j\pi_i$, for every $i,j$.

(a) We analyze the cases when $\pi_i=\alpha_l$ or $\beta_l$ and
$\pi_j=\alpha_h$ or $\beta_h$.

\emph{Case $(i)$.} If $\pi_i=\alpha_l$ and $\pi_j=\alpha_h$, then
$$\pi_i\pi_j=\sigma_l^+ A_{2l-1}A_{2l}\sigma_l^+  \sigma_h^+
A_{2h-1}A_{2h}\sigma_h^+ .$$ If $l=h$ then the claim is clearly
true. If $l \neq h$ we have that $\sigma_l^+
A_{2l-1}A_{2l}\sigma_l^+$ and $\sigma_h^+ A_{2h-1}A_{2h}\sigma_h^+
$ commute because they are disjoint permutations; hence the result
follows.

\emph{Case $(ii)$.} If $\pi_i=\beta_l$ and $\pi_j=\beta_h$. Idem.

\emph{Case $(iii)$.} If $\pi_i=\alpha_l$ and $\pi_j=\beta_h$, then
$$\pi_i\pi_j=\sigma_l^+ A_{2l-1}A_{2l}\sigma_l^+  \sigma_h^-
A_{2h-1}A_{2h}\sigma_h^-.$$ If $l \neq h$ then $\sigma_l^+
A_{2l-1}A_{2l}\sigma_l^+$ and $\sigma_h^-
A_{2h-1}A_{2h}\sigma_h^-$ commute; hence $\pi_i\pi_j=\pi_j\pi_i$.
While if $l=h$, it is easy to check that $\pi_i\pi_j= \alpha_l
\beta_l=\id=\beta_l\alpha_l=\pi_j\pi_i$, and the result follows.

(b) In general, for $\pi_i=(\sigma_{l_k}^{\pm}\cdots
\sigma_{l_1}^{\pm}) \trid \pi$ and
$\pi_j=(\sigma_{h_M}^{\pm}\cdots \sigma_{h_1}^{\pm}) \trid \pi$,
we have
\begin{align*}
\pi_i&=A_1A_2 \cdots \sigma_{l_1}^{\pm} A_{2l_1-1}A_{2l_1}
\sigma_{l_1}^{\pm} \cdots \sigma_{l_k}^{\pm} A_{2l_k-1}A_{2l_k}
\sigma_{l_k}^{\pm} \cdots A_{2L-1}A_{2L},\\
\pi_j&=A_1A_2 \cdots \sigma_{h_1}^{\pm} A_{2h_1-1}A_{2h_1}
\sigma_{h_1}^{\pm} \cdots \sigma_{h_M}^{\pm} A_{2h_M-1}A_{2h_M}
\sigma_{h_M}^{\pm} \cdots A_{2L-1}A_{2L}.
\end{align*}
We have two cases:\\
(i) If $l_r\neq h_s$, for all $r,s$. Then
\begin{align*}
\pi_i\pi_j&= \sigma_{l_1}^{\pm}
A_{2l_1-1}A_{2l_1}\sigma_{l_1}^{\pm}\dots \sigma_{l_k}^{\pm}
A_{2l_k-1}A_{2l_k}\sigma_{l_k}^{\pm}\\ & \quad  \sigma_{h_1}^{\pm}
A_{2h_1-1}A_{2h_1}\sigma_{h_1}^{\pm}\dots \sigma_{h_M}^{\pm}
A_{2h_M-1}A_{2h_M}\sigma_{h_M}^{\pm}=\pi_j\pi_i,
\end{align*}
because every $\sigma_{l_r}^{\pm}
A_{2l_r-1}A_{2l_r}\sigma_{l_r}^{\pm}$ commute with every
$\sigma_{h_s}^{\pm} A_{2h_s-1}A_{2h_s}\sigma_{h_s}^{\pm}$.\\
(ii) If $l_r = h_s$, for some $r,s$, we use (a) in every factor
corresponding to $l_r=h_s$. \epf

\smallbreak

For the rest of this subsection we fix the order in $T$ given by
$$
T=\{\pi_0=\pi,\pi_1=\alpha_1,\pi_2=\beta_1,\dots
,\pi_{2L-1}=\alpha_L,\pi_{2L}=\beta_L,\dots \}.
$$
Next we deal with $\rho=\rho_{\chi_{(k)},\mu}$ in
$\widehat{\snn^{\pi}}$, as in Section \ref{irrep}; as usual, let
$V$ be the vector space affording $\rho$ and $V_{\mu}$ the vector
space affording $\mu$. By Remark \ref{claim2bis}, we only need to
consider $k$ odd; thus $\rho(\pi)=-\Id$. Since $\pi_i^2=\id$, for
all $i$, then the possibles eigenvalues of the operators $\{
\rho(\pi_i) \, : \, 0 \leq i \leq N  \}$ are $1$ and $-1$.
Moreover, since $T$ is abelian  there exists a basis $\B$ of $V$
of simultaneous eigenvectors-- say  $\B = \{v_1,\dots,v_R\}$. Note
that $$\dim V = [\sn:\sn^{\chi_{(k)}}] \dim V_{\mu}=\binom{n}{k}
\dim V_{\mu}.$$

For every $i$, $0 \leq i \leq N$ we define ${\bf{f}}
^i=(f^i_1,f^i_2, \dots,f^i_R)$ where $\rho (\pi_i)v_r=f^i_r v_r$,
$1 \leq r \leq R$; for instance ${\bf{f}} ^0=(-1,-1, \dots,-1)$.
Now we denote by $E_i$ the matrix with all its rows equal to
${\bf{f}} ^i$. Hence $E_0$ is the matrix $\dim V \times \dim V$
with all its entries equal to $-1$.

Let us consider the subspace $W$ of $M(\Oc_{\pi}, \rho)$ with
basis $\{w_{i,r}:= \sigma_i v_r =\sigma_i \otimes v_r  \, : \,
0\leq i\leq N, 1\leq r \leq R \}$. Then $W$ is a braided vector
subspace of Cartan type and the matrix of the scalars
$(q_{a,b})_{a,b}$-- see section \ref{conventionsnichols}-- has the
form
$$
\mathcal Q=\begin{pmatrix} E_0 & E_1 & E_2 & \cdots &
E_{2L-1} &  E_{2L} & \cdots \\
E_1 & E_0 & E_1 & \cdots &
\cdots &  \cdots & \cdots \\
E_2 & E_2 & E_0 & \cdots &
\cdots &  \cdots & \cdots \\
\vdots & \vdots & \vdots & \ddots &
\cdots &  \cdots & \cdots \\
E_{2L-1} & \cdots & \cdots & \cdots &
E_0 &  E_{2L-1} & \cdots \\
E_{2L} & \cdots & \cdots & \cdots &
E_{2L} & E_0 & \cdots \\
\vdots & \vdots & \vdots & \vdots &
\vdots &  \vdots & \ddots \\
\end{pmatrix}.
$$
Here the diagonal blocks are equal to the matrix $E_0$; whereas
the block in the position $i,j$ is the matrix $E_k$ where $\pi_i
\sigma_j=\sigma_j \pi_k$.

\smallbreak\textbf{Notation in case $n$ odd.} We suppose $n=2L+1$
and for every $l$, with $1\leq l \leq L$, we take $\sigma_l^+$,
$\sigma_l^-$, $\alpha_l$, $\beta_l$ and $T$ as for $n$ even. So
\begin{align*}
\pi&=A_1A_2\dots A_{2l-1}A_{2l} \dots A_{2L-1}A_{2L}A_{2L+1},\\
\alpha_l&=A_1A_2\dots (4l-3 \,\, 4l-1)(4l-2 \,\, 4l) \dots A_{2L-1}A_{2L}A_{2L+1}, \\
\beta_l&=A_1A_2\dots (4l-3 \,\, 4l)(4l-2 \,\, 4l-1) \dots
A_{2L-1}A_{2L}A_{2L+1}.
\end{align*}
Then $\sigma_l^{\pm}$,$\alpha_l$, $\beta_l$, $\pi_i$, $\sigma_j$
and $T$ satisfy the same properties as before.

\begin{prop}\label{prop2}
Let $\rho=(\rho,V)$ be in $\widehat{\snn^{\pi}}$, $n \geq 2$. If
\begin{itemize}
\item[(a)] $n$ is even, or
\item[(b)] $n=3$ and $\rho \neq \chi_{-,\pm}$, or
\item[(c)] $n$ is odd and $\rho \neq \rho_{\chi_{(n)},\mu}$, for any
$\mu$ in $\widehat{\sn}$,
\end{itemize}
then $\dim\toba(\Oc_\trasp, \rho) = \infty$.
\end{prop}

\pf Let $\rho=\rho_{\chi_{(k)},\mu}$, $\chi_{(k)}$ in
$\widehat{\Z_2^n}$ and $\mu$ in $\widehat{\sn}$. By Remark
\ref{claim2bis}, we can assume that $k$ is odd. We distinguish two
possibilities.

(1) For every $l$, with $1\leq l\leq L$, $\rho(\alpha_l) =\Id
=\rho(\beta_l)$ or $\rho(\alpha_l) = -\Id = \rho(\beta_l)$. By
Lemma \ref{prop1} (1), this implies
$$\rho(A_{2l-1}A_{2l})=\Id, \qquad \text{ for all $l$.}$$

\emph{Assume that $n$ is even}. Hence $\rho(\pi)=\rho(A_1A_2
\cdots A_{2L-1}A_{2L})=\Id$; so $q_{\pi,\pi}=1$ and
$\dim\toba(\Oc_\trasp, \rho) = \infty$.

\emph{Assume that $n$ is odd}. Since $\rho(\pi)=-\Id$, it is easy
to see that $\rho(A_{2L+1})=-\Id$. By the discussion in subsection
\ref{irrep}, this implies $\rho(A_j)=-\Id$, $1 \leq j \leq 2L+1$.
Then $\rho=\rho_{\chi_{(n)},\mu}$, but this is a contradiction by
hypothesis.

\smallbreak (2) There exists $l$ with $\rho(\alpha_l)\neq \pm \Id$
or $\rho(\beta_l)\neq \pm \Id$ or $\rho(\alpha_l)= \pm \Id=\mp
\rho(\beta_l)$. Here we have that if $\dim V>4$ then the
generalized Cartan matrix $\mathcal A$ is such that its associated
Dynkin diagram is not of finite type and the result follows. For
see this we suppose that there exists $l$ with $\rho(\alpha_l)\neq
\pm \Id$; for the other cases the argument is similar. We regard
that the components of the vector ${\bf{f}}^l$ are $1$ or $-1$; so
we define $c^+:=\card \{ r : f^l_r=1 \}$ and $c^-:=\card \{ r :
f^l_r=-1 \}$; note that $c^+ +
c^-=R$. We consider three cases.\\
(i) If $\dim V \geq 7$ then the associated Dynkin diagram has a
vertex $w$ with $\lambda_w \geq 4$, where $\lambda_w$ denotes the
number of vertices of the diagram which are adjacent to $w$.
Hence, such diagram is not of finite type.\\
(ii) Let $\dim V = 6$; if $c^+ \geq 4$ or $c^- \geq 4$ we proceed
as in (i). So, we must consider $c^+ \leq 3$ and $c^- \leq 3$.
Because $c^+ + c^-=6$ then $c^+=3$ and $c^-=3$, but since there is
no Dynkin diagram of finite type with two vertices $w$, $w'$ with
$\lambda_w=3$ and $\lambda_{w'}=3$, the result follows.\\
(iii) If $\dim V = 5$ we only must consider either $c^+=3$ and
$c^-=2$ or $c^+=2$ and $c^-=3$, by (ii). In any case we have that
the associated Dynkin diagram has two vertices $w$, $w'$ with
$\lambda_w=3$ and $\lambda_{w'}=3$ and the result follows.
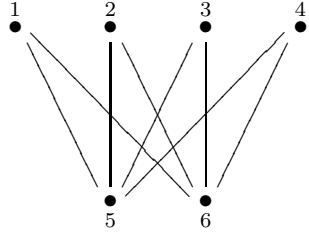
\[
\xymatrix{ {\overset{1}\bullet}\ar@{-}[2,1] \ar@{-}[2,2] &
{\overset{2}\bullet} \ar@{-}[2,0] \ar@{-}[2,1]
&{\overset{3}\bullet}\ar@{-}[2,-1]
\ar@{-}[2,0] &{\overset{4}\bullet} \ar@{-}[2,-1] \ar@{-}[2,-2] \\
\\&{\underset{5}\bullet} & {\underset{6}\bullet}& }
\]

Thus, we must consider $(\rho,V)$ with $\dim V \leq 4$. Then,
since $\dim V=\binom{n}{k} \dim V_{\mu}$, where $V_{\mu}$ is the
vector space affording $\mu$, we must consider the different
possibilities for $n$, $k$ and $\mu$ which satisfy the condition
$$\binom{n}{k} \dim V_{\mu}\leq 4.$$
This inequality holds only in the following cases
\begin{itemize}
\item[(i)] $n=2$ and $k=1$.
\item[(ii)] $n=3$, $k=1$ or $2$ and hence $\dim V_{\mu}=1$.
\item[(iii)] $n=4$, $k=1$ or $3$ and  hence $\dim V_{\mu}=1$.
\item[(iv)] any $n$, $k=0$ or $k=n$ and $\dim V_{\mu}=1,2,3$ or $4$.
\end{itemize}
In (i), (ii) and (iii) the result follows from \cite[Th. 2.7]{AZ},
Propositions \ref{dosdosdos} and \ref{dos4}, respectively. In the
case (iv), $k\neq 0$ by Remark \ref{claim2bis} and $k= n$ would be
considered for $n$ odd, but this was discarded by hypothesis. \epf

%%%%%%%%%%%%%%%%%%%%%%%%%%%%%%%%%%%%%%%%%%%%%%%%%%%%%%%%%%%%%%%%%
%%%%%%%%%%%%%%%%%%%%%%%%%%%%%%%%%%%%%%%%%%%%%%%%%%%%%%%%%%%%%%%%%

\begin{theorem}\label{dosn} Let $\trasp\in \snn$ of type $(2^{n})$.

(a). If $n$ is even then $\dim\toba(\Oc_\trasp, \rho) = \infty$
for any $\rho \in \widehat{\snn^{\trasp}}$.

(b). If $n$ is odd and $\rho\neq \chi_{(n)} \otimes \epsilon$,
$\chi_{(n)} \otimes \sgn$, then $\dim\toba(\Oc_\trasp, \rho) =
\infty$ for any $\rho \in \widehat{\snn^{\trasp}}$.
\end{theorem}

The braided vector spaces associated to $\chi_{(n)} \otimes
\epsilon$ or to $\chi_{(n)} \otimes \sgn$ were considered in
Proposition \ref{claim2}.

\pf We can assume that $\pi=(1\,2)(3\,4)\dots (2n - 1\, 2n)$. By
Propositions \ref{dosdosdos} and \ref{prop2}, we only need to
consider $3 < n$ odd  and $\rho=\rho_{\chi_{(n)},\mu}$, with $\mu$
in $\widehat{\sn}$, $\mu\neq\varepsilon, \sgn$. Notice that
$$\rho= \Ind_{\Z_2^{n} \rtimes \mathbb S_n^{\chi_{(n)}}}^{\Z_2^{n} \rtimes
\mathbb S_n} (\chi_{(n)}\otimes \mu)=\Ind_{\Z_2^{n} \rtimes
\mathbb S_n}^{\Z_2^{n} \rtimes \mathbb S_n} (\chi_{(n)}\otimes
\mu)=\chi_{(n)}\otimes \mu.$$ We distinguish two possibilities, as
in the proof of \ref{prop2}.

(1). We suppose $\rho(\alpha_l) =\Id =\rho(\beta_l)$ or
$\rho(\alpha_l) =-\Id =\rho(\beta_l)$, for every $l$, with $1\leq
l\leq L$. Then it is easy to check that $\rho(B_{2l-1})= \pm\Id $,
$1 \leq l \leq L$. Since $B_{2}, B_{4},\dots,B_{2L}$ are disjoint
permutations we have that the operators $\rho(B_j)$, $1 \leq j
\leq n$, commute. Hence, there exists a basis of simultaneous
eigenvectors of such operators. This says that the representation
$\mu$ is not irreducible unless $\dim V_{\mu}=1$, and therefore
$\mu= \sgn$, but this is a contradiction by hypothesis. The case
$\rho(\alpha_l) = -\Id = \rho(\beta_l)$, $1\leq l\leq L$, implies
$\rho(B_{2l-1})= \Id $, $1 \leq l \leq L$; by analogous arguments
we conclude $\mu= \epsilon$, a contradiction by hypothesis.

(2). If $n \geq 7$ and $\mu \neq \epsilon, \sgn$, then $\dim
V_{\mu} > 4$, see \cite[4.14]{FH}; hence
$\dim\toba(\Oc^{2n}_\trasp, \rho) = \infty$. It remains the case
$n=5$ and the representations
\begin{align*}
\rho=\chi_{(5)} \otimes \phi ,\qquad\rho=\chi_{(5)} \otimes \psi ,
\end{align*}
where $\phi$, $\psi$ are the two irreducible representations of
$\mathbb{S}_5$ of dimension 4, let us say $\phi$ the standard
representation of $\mathbb{S}_5$ and $\psi$ its conjugated
representation.

Let us consider $\rho=\chi_{(5)} \otimes \phi$. We take
$\pi=A_1A_2A_3A_4A_5$, $B_j$, $\sigma^{\pm}$, $\alpha_l$,
$\beta_l$, $\pi_i$, $\sigma_j$ and $T$ as in the case $n$ odd; so
$ \mathbb S_{10}^{\pi}=<A_1,A_2,A_3,A_4,A_5> \rtimes
<B_1,B_2,B_3,B_4> \simeq \Z_2^{5}\rtimes \mathbb S_5$ and $T=\{
\pi_0=\pi, \pi_1, \dots, \pi_8\}$ satisfying
$$
\pi_1=B_1A_3A_4A_5, \,\, \pi_2= \pi B_1, \,\, \pi_3=A_1A_2B_3A_5,
\,\, \pi_4= \pi B_3, \,\,   \pi_5=B_1B_3A_5,$$
$$ \pi_6=B_1A_3A_4B_3A_5, \,\, \pi_7=A_1A_2B_1B_3A_5, \,\, \pi_8=A_1A_2B_1
A_3A_4B_3A_5.
$$
It is easy to check that the standard representation of $\mathbb
S_5$ can be given by
$$
\phi(1 \, 2)=\small{\begin{pmatrix} -1 & -1 & -1 & -1\\ 0 & 1 & 0
& 0\\0 & 0 & 1 & 0 \\0 & 0 & 0 & 1
\end{pmatrix}}, \quad
\phi(2 \, 3)=\small{\begin{pmatrix} 0 & 1 & 0 & 0\\ 1 & 0 & 0 &
0\\0 & 0 & 1 & 0 \\0 & 0 & 0 & 1
\end{pmatrix}},
$$
$$
\phi(3 \, 4)= \small{\begin{pmatrix} 1 & 0 & 0 & 0\\ 0 & 0 & 1 &
0\\0 & 1 & 0 & 0 \\0 & 0 & 0 & 1
\end{pmatrix}}, \quad
\phi(4 \, 5)= \small{\begin{pmatrix} 1 & 0 & 0 & 0\\ 0 & 1 & 0 &
0\\0 & 0 & 0 & 1 \\0 & 0 & 1 & 0
\end{pmatrix}}.
$$
Then it is clear that
$$
\rho(\pi)=- \Id, \quad \rho(\pi_1)=\rho(\pi_2)=- \phi(B_1), \quad
\rho(\pi_3)=\rho(\pi_4)=- \phi(B_3),
$$
$$
\rho(\pi_5)=\rho(\pi_6)=\rho(\pi_7)=\rho(\pi_8)=-\phi(B_1)\phi(B_3)=\small{\begin{pmatrix}
1 & 1 & 1 & 1\\ 0 & 0 & -1 & 0\\0 & -1 & 0 & 0 \\0 & 0 & 0 & -1
\end{pmatrix}}.
$$
If $v_1=(1,0,0,0)$, $v_2=(0,1,-1,0)$, $v_3=(1,0,0,-2)$ and
$v_4=(1,1,1,-4)$ then they are simultaneous eigenvectors of those
operators. Hence, we have ${\bf{f}}^1={\bf{f}}^2=(1,-1,-1,-1)$,
${\bf{f}}^3={\bf{f}}^4=(-1,1,-1,-1)$ and
${\bf{f}}^5={\bf{f}}^6={\bf{f}}^7={\bf{f}}^8=(1,1,-1,-1)$. Thus,
in the basis
\begin{align*}
w_1&=\sigma_0 v_1, \quad &w_2&=\sigma_0 v_2, \quad &w_3&=\sigma_1
v_1\\
w_4&=\sigma_1 v_2, \quad &w_5&=\sigma_2 v_1, \quad &w_6&=\sigma_2
v_2,
\end{align*}
the braiding is diagonal of Cartan type and the matrix $\mathcal
Q$ of the scalars $\left( q_{a,b} \right) _{a,b}$ is
$$\mathcal Q =\begin{pmatrix}
-1 & -1 & 1 & -1 & 1 & -1\\ -1 & -1 & 1 & -1 & 1 & -1\\
1 & -1 & -1 & -1 & 1 & -1\\ 1 & -1 & -1 & -1 & 1 & -1\\
1 & -1 & 1 & -1 & -1 & -1 \\1 & -1 & 1 & -1 & -1 & -1
\end{pmatrix};
$$
the corresponding Dynkin diagram is $A^{(1)}_5$, which is not of
finite type. Hence $\dim\toba(\Oc^{10}_\trasp, \rho) = \infty$.

Finally, if $\rho=\chi_{(5)} \otimes \psi$ we proceed as the
previous case using that the representation $\psi$ is given by
$\psi= \sgn \times \phi$. So, in the same basis as before we have
that the braiding is diagonal of Cartan type and we obtain the
same matrix $\mathcal Q$; hence the result follows. \epf

\smallbreak

\subsection{Nichols algebras corresponding to even unmixed permutations}

Let $r,n \in \N$, $r,n \geq 2$. Let $\pi=A_1\dots A_n$ in $\mathbb
S_{2rn}$, where $A_j$ is the $2r$-cycle
$$
A_j=\left( 2rj-2r+1  \quad 2rj-2r+2 \quad \cdots \quad 2rj\right),
$$
for every $j$, $1\leq j \leq n$. As explained in section
\ref{conventionss3}, we have
\begin{equation}\label{S2rnpi}
\mathbb S_{2rn}^{\pi}=\langle A_1, \dots, A_n \rangle \rtimes
\langle B_1,\dots,B_{n-1}\rangle \simeq \Z_{2r}^n \rtimes \mathbb
S_n,
\end{equation}
where $B_i$ is the involution
$$
B_i= \left(2r(i-1)+1 \,\,\, 2ri+1 \right) \left(2r(i-1)+2 \,\,\,
2ri+2 \right) \cdots \left(2ri \,\,\, 2r(i+1) \right),
$$
$1\leq i \leq n-1$. Then $A_j$ and $B_i$ satisfy the relations
analogous to those in subsection \ref{NAntransp}. Let $\rho$ be an
irreducible representation of $\mathbb S_{2rn}^{\pi}$ of the form
\begin{equation}\label{formrho2}
\rho=\Ind_{\Z_{2r}^n \rtimes \mathbb S_n^{\chi}}^{\Z_{2r}^{n}
\rtimes \mathbb S_n} (\chi \otimes \mu),
\end{equation} where
$\chi \in \widehat{\Z_{2r}^n}$ and $\mu \in \widehat{\mathbb
S_n^{\chi}}$. Let $\omega=\exp(\frac{i \pi}{r}) \in \G_{2r}$ a
primitive $2r$-th root of $1$; any irreducible representation of
$\Z_{2r}^n$ is isomorphic to $\chi_{u_1,\dots,u_n}$, where
\begin{equation}\label{ecuacion2r}
 \chi_{u_1,\dots,u_n}(A_j)=\omega^{u_j}, \quad
\text{$1\leq j \leq n$},
\end{equation}
with $0\leq u_j \leq 2r-1$.

{\bf{Notation:}} if $\rho$ is as in \eqref{formrho2}, with $\chi$
as in \eqref{ecuacion2r}, we shall write
$\rho=\rho_{u_1,\dots,u_n,\mu}$.

By Lemma \ref{odd}, if $\rho(\pi) \neq -\Id$, then $\dim\toba(\C,
\rho)= \infty$. Hence, in the following we only consider
$\rho=\rho_{\chi_{u_1,\dots,u_n},\mu}$ such that $\rho(\pi)=
-\Id$; that is
\begin{align}\label{elomega}
\omega^{u_1+\cdots +u_n}=-1,
\end{align}
i.e. $u_1+\cdots +u_n=r$, $3r$, $5r$,\dots, $(2n-1)r$.

\smallbreak For every $(i,j)$, with $1 \leq i<j \leq n$, we define
$$
B_{ij}:=\begin{cases} B_i & \text{, if $|i-j|=1$},\\
B_i B_{i+1}\cdots B_{j-1}\cdots B_{i+1} B_i & \text{, if
$|i-j|>1$,}
\end{cases}
$$
and $\pi_{(i,j)}:=\pi B_{ij}$. Note that $B_{ij}$ acts as the
transposition $(i \,\,j)$ on $A_1,\dots,A_n$ and that
$\pi_{(i,j)}$ is in $\mathbb S_{2rn}^{\pi}$. We can state the
following.

\begin{lema} \label{lemaPRI}
For every $(i,j)$, with $1 \leq i<j \leq n$, we have
\begin{itemize}
\item[(a)] $\pi_{(i,j)}$ is in $\C$.
\item[(b)] there exists an involution $\sigma_{(i,j)}$ such that $\pi_{(i,j)}=\sigma_{(i,j)} \pi \sigma_{(i,j)}$.
\item[(c)] there exist involutions $\sigma$, $\widetilde{\sigma}_{(i,j)}$ in $\mathbb S_{2rn}$ such
that $\pi^{-1}=\sigma \pi \sigma$ and
$\pi_{(i,j)}^{-1}=\widetilde{\sigma}_{(i,j)} \pi
\widetilde{\sigma}_{(i,j)}$.
\end{itemize}
\end{lema}
\pf It is enough to prove this for $i=1$ and $j=2$.

(a) and (b). It is easy to see that
\begin{align*}
\pi_{(1,2)}=\pi &B_1= (1 \,\,\,\, 2r+2 \,\,\,\, 3 \,\,\,\, 2r+4
\,\,\,\, 5 \dots 4r-2 \,\,\,\, 2r-1 \,\,\,\, 4r) \\
&\times (2 \,\,\,\, 2r+3 \,\,\,\, 4 \,\,\,\, 2r+5 \,\,\,\, 6 \dots
4r-1 \,\,\,\, 2r \,\,\,\, 2r+1) A_3\cdots A_n,
\end{align*}
and we can choose
\begin{align*}
\sigma_{(1,2)}=(2 \,\,\,\, 2r+2 )(4 \,\,\,\, 2r+4 )\cdots (2r-2
\,\,\,\, 4r-2 )(2r \,\,\,\, 4r).
\end{align*}
Clearly, $\sigma_{(1,2)}$ is an involution and
$\pi_{(1,2)}=\sigma_{(1,2)} \pi \sigma_{(1,2)}\in \C$.

(c) For every $j$, $1 \leq j \leq n$,
$$ \sigma_{j}:=
\prod_{h=1}^{r}  \left( 2(j-1)r+h \,\,\,\, 2jr-h+1 \right)
$$
is an involution and satisfies $\sigma_{j}  \pi \sigma_{j}=A_1
\cdots A_j^{-1} \cdots A_n$. Now, if $\sigma=\sigma_{1} \cdots
\sigma_{n}$ then $\pi^{-1}=\sigma \pi \sigma$. Finally, \emph{if
$r$ is even} and
\begin{align*}
\widetilde{\sigma}_{(1,2)}=&(2 \,\,\,\, 4r)(4 \,\,\,\, 4r-2)\cdots
(2r \,\,\,\, 2r+2)\\ &(3 \,\,\,\, 2r-1)(5 \,\,\,\, 2r-3) \cdots
(r-3 \,\,\,\, r+5)(r-1
\,\,\,\, r+3)\\
&(2r+3 \,\,\,\, 4r-1)(2r+5 \,\,\,\, 4r-3)\cdots
(2r+r-1\,\,\,\,2r+r+3),
\end{align*}
or \emph{if $r$ is odd} and
\begin{align*}
\widetilde{\sigma}_{(1,2)}=&(2 \,\,\,\, 4r)(4 \,\,\,\, 4r-2)\cdots
(2r \,\,\,\, 2r+2)\\&(3 \,\,\,\, 2r-1)(5 \,\,\,\, 2r-3) \cdots
(r-2 \,\,\,\, r+4)(r \,\,\,\, r+2)\\ &(2r+3 \,\,\,\, 4r-1)(2r+5
\,\,\,\, 4r-3)\cdots (2r+r\,\,\,\,2r+r+2),
\end{align*}
then $\widetilde{\sigma}_{(1,2)}^2=\id$, and straightforward
computations imply that
$\pi_{(1,2)}^{-1}=\widetilde{\sigma}_{(1,2)} \pi
\widetilde{\sigma}_{(1,2)}$. \epf

\smallbreak We now consider two different cases according to the
degree of $\rho$.

\subsubsection{The degree of $\rho$ is greater that 1}\label{degrho>1}

\begin{theorem}\label{teordegrho>1}
Let $\rho$ be in $\widehat{\mathbb S_{2rn}^{\pi}}$. If
$\deg\rho>1$ then $\dim \toba(\C, \rho)= \infty$.
\end{theorem}

\pf

Let us consider two possibilities.

(A) \emph{Assume that there exists $(i,j)$, with $1 \leq i< j\leq
n$, such that $\rho\left( \pi_{(i,j)} \right)\neq \pm \Id$}. For
simplicity, we denote
\begin{align*}
t_1&:=\pi, \quad &t_2&:=\pi^{-1}, \quad &t_3&:=\pi_{(i,j)}, \quad
&t_4&:=\pi_{(i,j)}^{-1},\\
g_1&:=\id, \quad &g_2&:=\sigma, \quad &g_3&:=\sigma_{(i,j)}, \quad
&g_4&:=\widetilde{\sigma}_{(i,j)}.
\end{align*}
Now, we have the following relations: $t_1g_l=g_lt_l$ ,
$l=1,2,3,4$, and
\begin{align*}
t_2g_1&=g_1 t_2, \quad & t_2g_2&=g_2 t_1, \quad &t_2g_3&=g_3t_4,
\quad &t_2g_4&=g_4 t_3,\\
t_3g_1&=g_1t_3, \quad & t_3g_2&=g_2t_4, \quad & t_3g_3&=g_3t_1,
\quad & t_3g_4&=g_4t_2,\\
t_4g_1&=g_1t_4, \quad & t_4g_2&=g_2t_3, \quad & t_4g_3&=g_3t_2,
\quad &t_4g_4&=g_4t_1.
\end{align*}

Since the elements $t_1$, $t_2$, $t_3$ and $t_4$ commute then
there exists a basis of simultaneous eigenvectors $\{v_1,\dots,v_R
\}$ of $V$, the vector space affording $\rho$. Hence, either the
operator $\rho\left( \pi_{(i,j)} \right)$ has at least two
distinct eigenvalues or $\rho\left( \pi_{(i,j)} \right)=\lambda
\Id$, with $\lambda \neq \pm 1$.

In the first case, there exist $s$ and $s'$, $1 \leq s, s' \leq
R$, such that
$$\rho(\pi_{(i,j)}) \, v_s=\lambda_sv_s \quad \text{ and } \quad
\rho(\pi_{(i,j)}) \, v_{s'}=\lambda_{s'}v_{s'},$$ with $\lambda_s
\neq \lambda_{s'}$; let us consider the subspace $W$ of
$M(\C,\rho)$ generated by
\begin{equation}\label{subspaceW}
\{g_1 v_s, \, g_1 v_{s'},\, g_2 v_s, \, g_2 v_{s'}, \, g_3 v_s, \,
g_3 v_{s'}, \, g_4 v_s, \, g_4 v_{s'} \}.
\end{equation}
It is clear that $W$ is a braided vector subspace of diagonal type
of $M(\C,\rho)$. Now, if $\lambda_s^2\neq 1$ then it is easy to
see that the generalized Dynkin diagram contains a cycle of the
form
\begin{align}\label{4ciclo}
\xymatrix @+1.5pc @ur{-1 \bullet \ar@{-}[r]^{\lambda_s^2} &
{\overset{-1} \bullet} \ar@{-}[d]^{\lambda_s^{-2}}\\
{\underset{-1} \bullet} ´\ar@{-}[u]^{\lambda_s^{-2}} & \bullet \!
- \! 1 \ar@{-}[l]^{\lambda_s^2}};
\end{align}
while that if $\lambda_s^2=1$ then $\lambda_s\lambda_{s'} \neq 1$,
this implies that the generalized Dynkin diagram contains a cycle
of the form
\begin{align*}
\xymatrix @+1.5pc @ur{-1 \bullet
\ar@{-}[r]^{\lambda_s\lambda_{s'}} &
{\overset{-1} \bullet} \ar@{-}[d]^{\lambda_s^{-1}\lambda_{s'}^{-1}}\\
{\underset{-1} \bullet}
´\ar@{-}[u]^{\lambda_s^{-1}\lambda_{s'}^{-1}} & \bullet \! - \! 1
\ar@{-}[l]^{\lambda_s\lambda_{s'}}}.
\end{align*}
Hence, in both cases we have $\dim \toba(\C, \rho)= \infty$, by
\cite{H4}.

In the second case, we choose any $s$, $1\leq s \leq R$; then the
subspace of $M(\C,\rho)$ generated by
\begin{equation*}
\{g_1 v_s, \, \, g_2 v_s, \, \, g_3 v_s, \,
 \, g_4 v_s\},
\end{equation*}
is a braided vector subspace of diagonal type of $M(\C,\rho)$, and
its Dynkin diagram contains a cycle as in \eqref{4ciclo}. Hence,
$\dim \toba(\C, \rho)= \infty$, by \cite{H4}.

(B) \emph{Assume that $\rho\left( \pi_{(i,j)} \right)= \pm \Id$,
for every $(i,j)$, with $1 \leq i< j\leq n$}. The relation
$\pi_{(1,2)}=\pi B_1$ gives $\rho(B_1)=\mp \Id$; the relations
$\pi_{(1,3)}=\pi B_1B_2B_1$ and $\rho(B_1)=\pm\Id$ imply that
$\rho(B_2)=\mp \Id$, and so on. Hence, the operators $\rho(A_1)$,
\dots, $\rho(A_n)$, $\rho(B_1)$, \dots, $\rho(B_{n-1})$ commute,
and  there exists a basis of simultaneous eigenvectors of $V$ for
those operators. Since $\deg\rho >1$, $\rho$ is not an irreducible
representation of $\mathbb S_{2rn}^{\pi}$, which is a
contradiction. \epf

\bigbreak

\subsubsection{The degree of $\rho$ is 1} Say $V=\ku$ - span of $v$. By \eqref{formrho2},
$ \deg \rho=[\sn:\mathbb S_{n}^{\chi}] \,\, \deg \mu$; thus
$\mathbb S_{n}^{\chi}=\sn$ and $\deg \mu =1$. This implies that
$\rho=\chi_{c,\dots,c} \otimes \mu$, for some $c$, with $0\leq c
\leq 2r-1$, and $\mu= \epsilon$ or $\sgn$. Note that if $c=0$ then
$\rho(\pi)=1$, which is a contradiction by hypothesis. So, we can
assume $c\neq 0$.

We begin by the following result.

\begin{prop}\label{propfinal1} Let $\rho=\chi_{c,\dots,c} \otimes \mu$, with $0< c \leq
2r-1$.
\begin{itemize}
\item[(a)] If $r$ is odd and $c \neq r$, then $\dim \toba(\C,
\rho)= \infty$.
\item[(b)] If $r$ is even and $c \neq \frac{r}{2}$, $r$, $\frac{3r}{2}$, then $\dim \toba(\C,
\rho)= \infty$.
\end{itemize}
\end{prop}
\pf Let
\begin{align*}
t_1&:=\pi, \quad &t_2&:=\pi^{-1}, \quad &t_3&:=A_1^{-1}A_2\cdots
A_n, \quad
&t_4&:=t_3^{-1},\\
g_1&:=\id, \quad &g_2&:=\sigma, \quad &g_3&:=\sigma_1, \quad
&g_4&:=\sigma_2 \cdots \sigma_n,
\end{align*}
where $\sigma_1$, \dots, $\sigma_n$ are as in the proof of Lemma
\ref{lemaPRI} (c). It is clear that they  satisfy the same
relations as in subsection \ref{degrho>1}. Then the subspace of
$M(\C,\rho)$ generated by $\{g_1 v, \, g_2 v, \, g_3 v, \, g_4
v\}$ is braided of diagonal type which matrix of coefficients
$(q_{ij})_{ij}$, see subsection \ref{conventionsnichols}, given by
$$\mathcal Q =\begin{pmatrix}
-1 & -1 & -\omega^{-2c} & -\omega^{2c}\\ -1 & -1 & -\omega^{2c} & -\omega^{-2c}\\
-\omega^{-2c} & -\omega^{2c} & -1 & -1 \\ -\omega^{2c} &
-\omega^{-2c} & -1 & -1
\end{pmatrix}.
$$
Since $c\leq 2r-1$, it is easy to see that $\omega^{4c}=1$ if and
only if $2c=r$, $2r$ or $3r$. Now, it is clear that if $r$ is odd
and $c\neq r$, or if $r$ is even and $c \neq \frac{r}{2}$, $r$,
$\frac{3r}{2}$, we have that $\omega^{4c}\neq 1$. This implies
that the generalized Dynkin diagram has a cycle as in
\eqref{4ciclo}. Hence, $\dim \toba(\C, \rho)= \infty$. \epf

In the remaining cases, the braiding is always negative.

\begin{theorem}\label{teobra}
Assume that $\rho(\pi)=-1$.\\
(a) If $r$ is odd and $\rho=\chi_{r,\dots,r}\otimes \mu$, with
$\mu=\epsilon$ or $\sgn$, then the braiding is negative.\\
(b) If $r$ is even and $\rho=\chi_{c,\dots,c}\otimes \mu$, with
$c=\frac{r}{2}$, $r$ or $\frac{3r}{2}$ and $\mu=\epsilon$ or
$\sgn$, then the braiding is negative.
\end{theorem}
Note that, for $\rho$ as in (a) or (b), $\rho(\pi)$ is not
necessarily equal to $-1$.

In order to prove this result, we need two lemmata. Let us
remember that $t_1 = \pi$, \dots, $t_{M}$ is a numeration of $\C$
and $g_l\in \mathbb S_{2rn}$ are such that $g_l \pi g_l^{-1} =
t_l$, for all $1\le l \le M$; we choose $g_1=\id$.

Let $t_l$ in $\C$, such that $\pi t_l=t_l \pi$, i.e. $t_l$ in
$\mathbb S_{2rn}^{\pi}$. We know that $\gamma_{l1}:= g_1^{-1} t_l
g_1=t_l$ and  $\gamma_{1l}:= g_l^{-1}\pi g_l$ are in $\C \cap
\mathbb S_{2rn}^{\pi}$. By \eqref{S2rnpi}, we can write
\begin{align}
\gamma_{l1}&\label{gamma1}=A_1^{d_1}\cdots A_n^{d_n} B,\\
\gamma_{1l}&\label{gamma2}=A_1^{e_1}\cdots A_n^{e_n} B',
\end{align}
where $B$ and $B'$ are in $\langle B_1,\dots,B_{n-1} \rangle\simeq
\sn$. Let $\Phi : \langle B_1,\dots,B_{n-1} \rangle \to \sn$ be
the group isomorphism given by $\Phi(B_i)=(i \, \, i+1)$, $1 \leq
i \leq n-1$.

For every $j$, $1\leq j \leq n$, we define
$$\mathbb A_j:=\{2rj-2r+1, 2rj-2r+2 , \dots , 2rj \} ,$$
i.e. $\mathbb A_j$ is the set of natural numbers that are ``moved"
by $A_j$. We also set
\begin{equation}\label{J}
J:=\{j \, | \, A_jB\neq BA_j  \}.
\end{equation}
If $j\not \in J$ then $d_j$ is relatively prime to $2r$, because
$A_j^{d_j}$ is a cycle of length $2r$. Clearly, if the type of
$\Phi(B)$ is $(L)$ then $\card J= L$. So, we can write
\begin{equation}\label{ecua}
\gamma_{l1}= \prod_{j \not\in J}A_j^{d_j} \quad \prod_{j \in
J}A_j^{d_j}B,
\end{equation}
and it is easy to see that $g_l$ can be chosen
\begin{equation}\label{gl}
g_l= \nu \, \, \prod_{j\not\in J} \sigma_{l,j} \,\, ,
\end{equation}
where $\sigma_{l,j}A_j\sigma_{l,j}^{-1}=A_j^{d_j}$ and every
element of $\mathbb A_{j'}$ is fixed by $\sigma_{l,j}$ if
$j\not\in J$ and $j' \neq j$, and $\nu$ is such that
\begin{equation}\label{nu}
\nu \, \, \prod_{j \in J}A_j \, \, \nu^{-1}=\prod_{j \in
J}A_j^{d_j}B,
\end{equation}
and that every element of $\mathbb A_{j}$, $j\not \in J$, is fixed
by $\nu$.

\begin{lema}\label{lemafin1}
$\Phi(B)$ and $\Phi(B')$ have the same type in $\sn$.
\end{lema}

\pf We will consider cases according to the type of $\Phi(B)$ in
$\sn$.

If the type of $\Phi(B)$ is $(1^{2rn})$; this means $B=\id$. We
have that $J=\emptyset$, so we can chose
$g_l=\sigma_{l,1}\cdots\sigma_{l,n}$. Then
$$
\gamma_{1l}=g_l^{-1} \pi g_l=\sigma_{l,1}^{-1} A_1\sigma_{l,1}
\cdots \sigma_{l,n}^{-1}A_n\sigma_{l,n},
$$
and since $\gamma_{1l}$ is in $\C$, i.e. it is a product of
disjoint cycles of length $2r$, we have that $\sigma_{l,j}^{-1}
A_j\sigma_{l,j}$ is a cycle of length $2r$, for all $j$. This
implies that
$$
\gamma_{1l}=A_1^{e_1}\cdots A_n^{e_n},
$$
with $e_1,\dots,e_n$ relatively primes to $2r$; this means that
$B'=\id$.

If the type of $\Phi(B)$ is $(2)$. It is enough to assume that
$B=B_i$ for some $i$, $1 \leq i \leq n-1$. We saw that if $j\neq
i, i+1$ then $d_j$ is relatively prime to $2r$, and that $g_l$ can
be chosen as in \eqref{gl}, i.e.
$$
g_l= \nu  \prod_{j\neq i,i+1} \sigma_{l,j} \,\, ,
$$
where $\nu$ satisfies $\nu A_iA_{i+1}
\nu^{-1}=A_i^{d_i}A_{i+1}^{d_{i+1}}B,$ and if $j\neq i,i+1$ then
the elements of $\mathbb A_j$ are fixed by $\nu$. Hence,
\begin{align*}
\gamma_{1l}&= g_l^{-1} \pi g_l=\prod_{j\neq i,i+1}
\sigma_{l,j}^{-1} A_j\sigma_{l,j}  \quad \nu^{-1}  A_iA_{i+1} \nu
=\prod_{j\neq i,i+1} A_j^{e_j} \quad A_i^{e_i}A_{i+1}^{e_{i+1}}
B',
\end{align*}
with $e_j$ relatively prime to $2r$, if $j \neq i,i+1$. This
implies that the type of $\Phi(B')$ is
$(h_1^{m_1},\dots,h_K^{m_K})$ with
$$
m_1 h_1+ \cdots+m_K h_K\leq 2.
$$
Then the type of $\Phi(B')$ is $(1)$ or $(2)$; if it is $(1)$ we
have that $B'=\id$, then $B=\id$, by the first case, a
contradiction. Thus, the type of $\Phi(B')$ is $(2)$.

\smallbreak

Notice that if the type of $\Phi(B)$ is $(2^a)$ then the same
occurs for $\Phi(B')$, by repeating the previous argument in each
disjoint transposition that appears in the decomposition of
$\Phi(B)$ as product of disjoint permutations of $\sn$.

In general, we can prove by the same argument that if the result
is true when the type of $\Phi(B)$ is $(L_1)$ and $(L_2)$ then the
result is also true if the type of $\Phi(B)$ is $(L_1^2)$, if
$L_1=L_2$, or $(L_1,L_2)$, if $L_1 \neq L_2$.

Let $\Phi(B)$ be of type $(L)$. We use induction in $L$ and the
previous paragraphs to prove that the type of $\Phi(B')$ is $(L)$.
Explicitly, assume that there exists $L>2$ such that for every
$h<L$ it is true the following: if the type of $\Phi(B)$ is $(h)$,
then the type of $\Phi(B')$ is $(h)$. Suppose that the type of
$\Phi(B')$ is $(h_1^{m_1},\dots,h_K^{m_K})$. We proceed as in the
case $L=2$. We can chose $g_l$ as in \eqref{gl}, with $\nu$ that
satisfies \eqref{nu}, and if $j\not\in J$ then the elements of
$\mathbb A_j$ are fixed by $\nu$. Hence
\begin{align*}
\gamma_{1l}&= g_l^{-1} \pi g_l=\prod_{j\not\in J}
\sigma_{l,j}^{-1} A_j\sigma_{l,j}  \quad \nu^{-1} \prod_{j\in J}
A_j \, \nu =\prod_{j\not\in J} A_j^{e_j} \quad \prod_{j\in J}
A_j^{e_j} B',
\end{align*}
with $e_j$ relatively prime to $2r$ if $j \not\in J$, because
$\gamma_{1l}$ is in $\C$. This implies that
$$
m_1 h_1+ \cdots+m_K h_K\leq L.
$$
If $m_1 h_1+ \cdots+m_K h_K <L$ or if $m_1 h_1+ \cdots+m_K h_K =L$
with $K>1$, then $h_1,\dots,h_K<L$, and by inductive hypothesis
and the previous paragraph we have that the type of $\Phi(B)$ is
$(h_1^{m_1},\dots,h_K^{m_K})\neq (L)$, which is a contradiction.
So, type of $\Phi(B')$ is $(h_1)^{m_1}$, with $m_1 h_1=L$; if
$m_1>1$ we use inductive hypothesis and the previous paragraph to
say that the type of $\Phi(B)$ is $(h_1^{m_1})\neq (L)$, which is
a contradiction. Hence, $m_1=1$ and $h_1=L$, it means that the
type of $\Phi(B')=(L)$, and this concludes the proof. \epf

\begin{lema}\label{lemafin2}  Let $\gamma_{l1}$ and $\gamma_{1l}$ as in
\eqref{gamma1} and \eqref{gamma2}, respectively.\\
(a) For any $r$ if $n$ is odd, then $\sum_{j=1}^n (e_j + d_j)$ is even.\\
(b) If $r$ is even and $n$ is even, then $\sum_{j=1}^n (e_j +
d_j)\equiv 0 \mod(4)$.
\end{lema}
\pf (a) If $n$ is odd we have that the sign of $\pi$ in $\mathbb
S_{2rn}$ is
$$
\sgn \pi= \sgn A_1 \cdots \sgn A_n=(-1)^n=-1,
$$
because $A_1,\dots,A_n$ are cycles of even length. Since
$\gamma_{l1} \in \C$ we have that $\sgn \gamma_{l1}=-1$, on the
other hand
$$
\sgn \gamma_{l1}=\sgn A_1^{d_1} \cdots  \sgn A_n^{d_n} \sgn
B=(-1)^{d_1+\cdots+d_n},
$$
because $B \in \langle B_1,\dots,B_{n-1}\rangle$ and every
$B_1,\dots,B_{n-1}$ is a product of an even number of
transpositions in $\mathbb S_{2rn}$. Then $d_1+\cdots+d_n$ is odd.
Analogously, $e_1+\cdots+e_n$ is odd. Then the result follows.

(b) Assume that $n$ is even. In this case the sign of $\pi$ in
$\mathbb S_{2rn}$ is 1; since $\gamma_{l1}$ and $\gamma_{1l}$ are
in $\C$, $d_1+\cdots+d_n$ and $e_1+\cdots+e_n$ are even. We
suppose that the decomposition of $\Phi(B)$ in product of disjoint
permutation in $\sn$ is
\begin{equation}\label{desc}
\Phi(B)= \tau_1 \cdots \tau_K.
\end{equation}
By Lemma \ref{lemafin1}, we have that
\begin{equation*}
\Phi(B')= \tau_1' \cdots \tau_K',
\end{equation*}
with $|\tau_k|=|\tau_k'|$, for all $k$. Obviously, $|B|=$\,
lcm$(|\tau_1|,\dots,|\tau_K|)=|B'|$.

For every $k$, $1 \leq k \leq K$, we define
\begin{align}
J_k:=\{j  \, | \, 1\leq j \leq n \text{ and }
A_j\Phi^{-1}(\tau_k)\neq \Phi^{-1}(\tau_k)A_j  \}.
\end{align}
Clearly, $\card J_k=|\tau_k|$, for all $k$. Note that
$J_1,\dots,J_K$ are disjoint sets and if $J$ is as in \eqref{J}
then $J=J_1\cup \cdots \cup J_K$. Besides, it is clear that
$$
J_k=\{j  \, | \, 1\leq j \leq n \text{ and }
A_j\Phi^{-1}(\tau_k')\neq \Phi^{-1}(\tau_k')A_j  \}.
$$
by Lemma \ref{lemafin1}. We write $\gamma_{l1}$ as in \eqref{ecua}
in a more precise form
\begin{align*}
\gamma_{l1}=g_l \pi g_l^{-1}=\prod_{j\not\in J}A_j^{d_j} \quad
\prod_{j \in J_1}A_j^{d_j}\cdots \prod_{j \in J_K}A_j^{d_j} \quad
B,
\end{align*}
and $g_l$ can be chosen as in \eqref{gl}
\begin{align*}
g_l= \nu_1 \cdots \nu_K \,\, \prod_{j\not\in J} \sigma_{l,j}\,\, ,
\end{align*}
where
\begin{align}\label{paratruco1}
\nu_k \prod_{j\in J_k} A_j \,\, \nu_k^{-1}=\prod_{j\in
J_k}A_j^{d_j} \,\, \Phi^{-1}(\tau_k),
\end{align}
and if $j\not\in J_k$ every element of $\mathbb A_j$ is fixed by
$\nu_k$; this allows to say that if $j\not\in J_k$ then $A_j$ and
$\nu_k$ commute. Hence, if $\gamma_{1l}$ is as in \eqref{gamma2}
then
\begin{align*}
\gamma_{1l}=g_l \pi g_l^{-1}=\prod_{j\not\in J}A_j^{e_j} \quad
\prod_{j \in J_1}A_j^{e_j}\cdots \prod_{j \in J_K}A_j^{e_j} \quad
B',
\end{align*}
with
\begin{align}\label{paratruco2}
\prod_{j \in J_k}A_j^{e_j} \Phi^{-1}(\tau_k')= \nu_k^{-1} \prod_{j
\in J_k}A_j \,\, \nu_k.
\end{align}

Since $|\gamma_{l1}|=2r$, then $B^{2r}=\id$; this implies that
$|B|$ divides $2r$, let us say $2r=|B|q$, with $q\geq 1$. It is
straightforward to prove that
\begin{align}\label{importante}
\left( \prod_{j \in J}A_j^{d_j}B \right)^{h|B|}=\left( \prod_{j
\in J_1}A_j \right)^{h\frac{|B|}{|\tau_1|}\sum_{j\in J_1}d_j}
\cdots \left( \prod_{j \in J_K}A_j
\right)^{h\frac{|B|}{|\tau_K|}\sum_{j\in J_K}d_j},
\end{align}
for all integer $h$. When $h=q$ both sides are equal to $\id$ and
this implies
\begin{equation*}
\left( \prod_{j \in J_k}A_j
\right)^{q\frac{|B|}{|\tau_k|}\sum_{j\in J_k}d_j}=\id,
\end{equation*}
for all $k$. Since the order of $\prod_{j \in J_k}A_j$ is $2r$ we
have that $|\tau_k|$ divides $\sum_{j\in J_k}d_j$. Analogously, we
can prove that $|\tau_k'|$ divides $\sum_{j\in J_k}e_j$, for all
$K$. Hence, for every $k$, $1\leq k \leq K$, there exist
$p_k,p_k'\geq 1$ such that
\begin{align}\label{ecuap}
\sum_{j\in J_k}d_j= |\tau_k| p_k \qquad \text{ and } \qquad
\sum_{j\in J_k}e_j= |\tau_k| p_k'.
\end{align}

By \eqref{paratruco1}, \eqref{paratruco2} and \eqref{importante},
for every $k$ we have that
\begin{align*}
\left( \prod_{j \in J_k}A_j \right)^{h |B|}&=\nu_k^{-1} \nu_k
\left( \prod_{j \in J_k}A_j \right)^{h |B|}
\nu_k^{-1}\nu_k=\nu_k^{-1} \left( \nu_k  \prod_{j \in J_k}A_j
\nu_k^{-1}
\right)^{h |B|} \nu_k\\
&=\nu_k^{-1} \left( \prod_{j\in J_k}A_j^{d_j} \Phi^{-1}(\tau_k)
\right)^{h |B|} \nu_k
=\nu_k^{-1}\left( \prod_{j\in J_k}A_j  \right)^{h |B|p_k} \nu_k\\
&=\left( \nu_k^{-1} \prod_{j\in J_k}A_j   \nu_k \right)^{h
|B|p_k}=\left( \prod_{j \in J_k}A_j^{e_j} \Phi^{-1}(\tau_k')
\right)^{h
|B|p_k}\\
&=\left( \prod_{j\in J_k}A_j  \right)^{h |B|p_kp_k'},
\end{align*}
for all integer $h$. In particular, for $h=1$ this implies that
$2r$ divides $|B|p_kp_k'-|B|$. Since $2r=|B|q$, we have that $q$
divides $p_kp_k'-1$, for every $k$; let us say that for every $k$,
$1\leq k\leq K$, there exists $x_k\geq 1$ such that
\begin{align}\label{qdivide}
p_kp_k'-1= q x_k.
\end{align}

By a similar argument as in the previous paragraph, we can show
that
\begin{align}
\left( \prod_{j \in J_k}A_j \right)^{h |\tau_k|}=\left( \prod_{j
\in J_k}A_j \right)^{h |\tau_k|p_kp_k'},
\end{align}
for all integer $h$. For $h=1$, this says that $2r$ divides
$|\tau_k|p_kp_k'-|\tau_k|$. Using \eqref{qdivide} and that
$|B|=|\tau_k| y_k$, for some $y_k\geq 1$, we have that $y_k$
divides $x_k$, for every $k$, let us say $x_k=y_kz_k$, for some
$z_k\geq 1$. Hence
\begin{align*}
|\tau_k|p_kp_k'-|\tau_k|=|\tau_k| \, q \, y_k \, z_k=|B|\, q \,
z_k=2rz_k.
\end{align*}
Since $r$ is even we have that $|\tau_k|p_kp_k' \equiv |\tau_k|
\mod (4)$; this means that
\begin{align}\label{cuadrado1}
p_k' \sum_{j\in J_k}d_j  \equiv |\tau_k| \mod (4) \quad \text{and}
\quad p_k \sum_{j\in J_k}e_j  \equiv |\tau_k| \mod (4).
\end{align}
Clearly,
\begin{align}\label{cuadrado2}
(\sum_{j\in J_k}d_j) \,\, (\sum_{j\in J_k}e_j)  \equiv |\tau_k|^2
\mod (4).
\end{align}
Using \eqref{cuadrado1}, \eqref{cuadrado2} and that
$|\tau_k|^2\equiv 0 \text{ or } 1 \mod (4)$ we conclude that
\begin{align}
\sum_{j\in J_k}d_j \equiv \sum_{j\in J_k}e_j   \mod (4),
\end{align}
for every $k$, $1\leq k \leq K$. Moreover, if $|\tau_k|^2\equiv 0
\mod (4)$, then
\begin{align}\label{congruenciapar}
\sum_{j\in J_k}d_j \equiv  0 \equiv  \sum_{j\in J_k}e_j   \mod (4)
\quad \text{or} \quad \sum_{j\in J_k}d_j \equiv  2 \equiv
\sum_{j\in J_k}e_j \mod (4),
\end{align}
and if $|\tau_k|^2\equiv 1 \mod (4)$, then
\begin{align}\label{congruenciaimpar}
\sum_{j\in J_k}d_j \equiv  1 \equiv  \sum_{j\in J_k}e_j   \mod (4)
\quad \text{or} \quad \sum_{j\in J_k}d_j \equiv 3 \equiv
\sum_{j\in J_k}e_j \mod (4).
\end{align}

For $h=0$, $1$, $2$ and $3$, we define
\begin{align*}
\mathcal K_h:=\{k \, |  \, 1\leq k \leq K \text{ and } |\tau_k|^2
\equiv h \mod (4) \}.
\end{align*}
Then we can write
\begin{align*}
\sum_{j=1}^n(e_j+d_j)=\sum_{j\not\in J}(e_j+d_j)+\sum_{k\in
\mathcal K_1\cup \mathcal K_3}\sum_{j\in J_k}(e_j+d_j)+ \sum_{k\in
\mathcal K_0\cup \mathcal K_2}\sum_{j\in J_k}(e_j+d_j).
\end{align*}
By \eqref{congruenciapar}, it is clear that
$$
\sum_{k\in \mathcal K_0\cup \mathcal K_2}\sum_{j\in
J_k}(e_j+d_j)\equiv 0\mod (4),
$$
while if $k\in K_1\cup \mathcal K_3$ then $\sum_{j\in
J_k}(e_j+d_j)\equiv 2 \mod (4)$. Besides, if $j\not\in J$ then
$d_j$ and $e_j$ are relatively prime to $2r$ and it is easy to see
that
\begin{align}
A_j=\sigma_{l,j}^{-1}\sigma_{l,j}A_j\sigma_{l,j}^{-1}\sigma_{l,j}=\sigma_{l,j}^{-1}
A_j^{d_j}\sigma_{l,j}=(\sigma_{l,j}^{-1}
A_j\sigma_{l,j})^{d_j}=A_j^{e_jd_j};
\end{align}
this implies that $2r$ divides $e_jd_j-1$, and since $r$ is even
we have $e_jd_j \equiv 1\mod (4)$. Using that $d_j$ and $e_j$ are
odd and the last fact we can prove that
\begin{align}
e_j+d_j \equiv 2 \mod (4),
\end{align}
for every $j\not\in J$.

We saw that $\card J_k=|\tau_k|$, hence $\sum_{k\in K_0\cup
\mathcal K_2}\card J_k$ is even. Since
$$n=\card J^c+  \sum_{k\in K_1\cup \mathcal
K_3}\card J_k+\sum_{k\in K_0\cup \mathcal K_2}\card J_k$$ is even
we have that $a:=\card J^c+  \sum_{k\in K_1\cup \mathcal K_3}\card
J_k$ is even. Hence
\begin{align}
\sum_{j=1}^n(e_j+d_j) \equiv 2^a \equiv 0 \mod (4),
\end{align}
and the result follows. \epf

%-) Si no que da claro esto ultimo se puede hidratar escribiendo
% lo que esta en las hojas
%-) Ver que no genere problemas el caso $J=\emptyset$. Tratarlo
%aparte.

\emph{Proof of Theorem \ref{teobra}} Let $t_h$, $t_l$ in $\C$ that
commute; it amounts to say that $\gamma_{hl}:=g_l^{-1} t_h g_l$
and $\gamma_{lh}:=g_h^{-1} t_l g_h$ are in $\mathbb
S_{2rn}^{\pi}$. Let $\rho \in  \widehat{\mathbb S_{2rn}^{\pi}}$ as
in the statements (a) or (b). Let us remember from subsection
\ref{conventionsnichols}, that $q_{hh}=\rho(\gamma_{hh})$,
$q_{ll}=\rho(\gamma_{ll})$, $q_{hl}=\rho(\gamma_{hl})$ and
$q_{lh}=\rho(\gamma_{lh})$. We must see that $q_{hh}=-1=q_{ll}$
and $q_{hl}q_{lh}=1$. The first conditions are trivially
fulfilled. For the last one we consider two cases.

\emph{CASE 1:} $h=1$. Let $\gamma_{l1}$, $\gamma_{1l}$ be as in
\eqref{gamma1} and \eqref{gamma2}, respectively.

(i) \emph{Assume that $\rho=\chi_{r,\dots,r}\otimes \mu$, with
$\mu=\epsilon$ or $\sgn$}. Since $\rho(\pi)=-1$ and $\rho(\pi)=
\omega^{rn}$, with $\omega=\exp(\frac{i \pi}{r})$, we have that
$n$ must be odd. Then
\begin{align*}
q_{1l}q_{l1}&=(\chi_{r,\dots,r}\otimes \mu)
(\gamma_{1l}\gamma_{l1}) =(\chi_{r,\dots,r}\otimes \mu)
(\gamma_{1l})\,\,
(\chi_{r,\dots,r}\otimes \mu) (\gamma_{l1})\\
&=\omega^{r \sum_{j=1}^n e_j+d_j} \mu (B') \mu (B)=
(-1)^{\sum_{j=1}^n e_j+d_j}=1,
\end{align*}
by Lemma \ref{lemafin1} and Lemma \ref{lemafin2} (a).

(ii) \emph{Assume that $r$ is even and
$\rho=\chi_{c,\dots,c}\otimes \mu$, with $c=\frac{r}{2}$ or
$\frac{3r}{2}$, and $\mu=\epsilon$ or $\sgn$}. The condition
$\rho(\pi)=-1$ implies that $n \equiv 2 \mod (4)$; in particular
$n$ is even. By Lemma \ref{lemafin1} and Lemma \ref{lemafin2} (b),
we can say
\begin{align*}
q_{1l}q_{l1}&=(\chi_{c,\dots,c}\otimes \mu)
(\gamma_{1l}\gamma_{l1}) =(\chi_{c,\dots,c}\otimes \mu)
(\gamma_{1l})\,\,
(\chi_{c,\dots,c}\otimes \mu) (\gamma_{l1})\\
&=\omega^{c \sum_{j=1}^n e_j+d_j} \mu (B') \mu (B)= (\pm
i)^{\sum_{j=1}^n e_j+d_j}=1.
\end{align*}

\emph{CASE 2:} $h \neq 1$. We call $\tilde{\pi}:=t_h$.

(i) \emph{Assume that $\rho=\chi_{r,\dots,r}\otimes \mu$, with
$\mu=\epsilon$ or $\sgn$}. Then there exists $\tilde{c}$, with $0
\leq \tilde{c} \leq 2r-1$, such that
$\chi_{\tilde{c},\dots,\tilde{c}}\otimes \widetilde{\mu} \in
\widehat{\mathbb S _{2rn}^{\tilde{\pi}}}$ and
\begin{align}
M(\C,\chi_{r,\dots,r}\otimes \mu)=
M(\C,\chi_{\tilde{c},\dots,\tilde{c}}\otimes \widetilde{\mu}),
\end{align}
where $\widetilde{\mu}=\epsilon$ or $\sgn$, say
$\widetilde{\mu}=\sgn$. This implies that
$\widetilde{\rho}:=\chi_{\tilde{c},\dots,\tilde{c}}\otimes
\widetilde{\mu}$ and $\rho$ have the same image, see
\eqref{image}, it means that $\langle \pm
\omega^{\tilde{c}}\rangle=\langle \omega^r \rangle =\{1,-1 \}$.
Since $\omega=\exp(\frac{i \pi}{r})$ it is clear that
$\tilde{c}=r$. Now, the result follows for $\widetilde{\rho}$ from
the case (1)(i). The case $\widetilde{\mu}=\epsilon$ is similar.

(ii) \emph{Assume that $r$ is even and
$\rho=\chi_{c,\dots,c}\otimes \mu$, with $c=\frac{r}{2}$ or
$\frac{3r}{2}$, and $\mu=\epsilon$ or $\sgn$}. Then there exists
$\tilde{c}$, with $0 \leq \tilde{c} \leq 2r-1$, such that
\begin{align}
M(\C,\chi_{c,\dots,c}\otimes \mu)=
M(\C,\chi_{\tilde{c},\dots,\tilde{c}}\otimes \widetilde{\mu}),
\end{align}
where $\widetilde{\mu}=\epsilon$ or $\sgn$, say
$\widetilde{\mu}=\sgn$. This implies that
$\widetilde{\rho}:=\chi_{\tilde{c},\dots,\tilde{c}}\otimes
\widetilde{\mu}$ and $\rho$ have the same image, i.e. $\langle \pm
\omega^{\tilde{c}}\rangle=\langle \omega^r \rangle =\{1,i,-1,-i
\}$. Since $\omega=\exp(\frac{i \pi}{r})$ it is clear that
$\tilde{c}=\frac{r}{2}$ or $\frac{3r}{2}$. Now, the result follows
for $\widetilde{\rho}$ from the case (1)(ii). The case
$\widetilde{\mu}=\epsilon$ is similar.

This concludes the proof. \qed

%Now, we use that $\deg \rho=1$ and same argument as in the proof
%of Proposition \ref{claim2}

\end{document}